\documentclass[11pt]{article}

\usepackage{amsmath,latexsym,amssymb,amsfonts}

\newtheorem{theorem}{Theorem}[section]
\newtheorem{proposition}[theorem]{Proposition}
\newtheorem{corollary}[theorem]{Corollary}
\newtheorem{lemma}[theorem]{Lemma}
\newtheorem{remark}[theorem]{Remark}
\newtheorem{example}[theorem]{Example}

\newtheorem{notation}[theorem]{Notation}

\newenvironment{proof}[1][Proof.]{\begin{trivlist}
\item[\hskip \labelsep {\bfseries #1}]}{\end{trivlist}}

\def\bB{ {\mathbb B} }
\def\bC{ {\mathbb C} }
\def\cM{ {\mathcal M} }
\def\bR{ {\mathbb R} }

\def\id{{\mathbb{B}}}

\begin{document}

\title{On a remarkable semigroup of homomorphisms \\
with respect to free multiplicative convolution}

\author{Serban T. Belinschi 
\and Alexandru Nica  \thanks{Research supported by a Discovery Grant 
of NSERC, Canada and by a PREA award from the province of Ontario.}  }

\date{ }

\maketitle

\vspace{.5cm}

\begin{abstract}
Let $\cM$ denote the space of Borel probability measures on $\bR$.
For every $t \geq 0$ we consider the transformation 
$\bB_t : \cM \to \cM$ defined by 
\[
{\mathbb B}_t ( \mu ) = \Bigl( \, \mu^{\boxplus (1+t)} \, 
\Bigr)^{\uplus (1/(1+t))}, \ \ \mu \in {\mathcal M} ,
\]
where $\boxplus$ and $\uplus$ are the operations of free additive 
convolution and respectively of Boolean convolution on $\cM$, and where 
the convolution powers with respect to $\boxplus$ and $\uplus$ are defined 
in the natural way. We show that $\bB_s \circ \bB_t = \bB_{s+t}$, 
$\forall \, s,t \geq 0$ and that, quite surprisingly, every $\bB_t$ is a 
homomorphism for the operation of free {\em multiplicative} convolution 
$\boxtimes$ (that is, $\bB_t ( \mu \boxtimes \nu )$ =
$\bB_t ( \mu ) \boxtimes \bB_t ( \nu )$ for all $\mu , \nu \in \cM$ such that
at least one of $\mu , \nu$ is supported on $[0, \infty )$).

We prove that for $t=1$ the transformation $\bB_1$ coincides with 
the canonical bijection $\bB : \cM \to \cM_{\mathrm{inf-div}}$ discovered 
by Bercovici and Pata in their study of the relations between
infinite divisibility in free and in Boolean probability. 
Here $\cM_{\mathrm{inf-div}}$ stands for the set of probability distributions 
in $\cM$ which are infinitely divisible with respect to the operation $\boxplus$. 
As a consequence, we have that $\bB_t ( \mu )$ is $\boxplus$-infinitely divisible 
for every $\mu \in \cM$ and every $t \geq 1$.

On the other hand we put into evidence a relation between the transformations 
$\bB_t$ and the free Brownian motion; indeed, Theorem 4 of the paper gives an
interpretation of the transformations ${\mathbb B}_t$ as a way of re-casting the 
free Brownian motion, where the resulting process becomes multiplicative with 
respect to $\boxtimes$, and always reaches $\boxplus$-infinite divisibility by 
the time $t=1$.
\end{abstract}

\newpage

\section{Introduction}

$\ $

{\bf 1.1 The transformations \boldmath{${\mathbb B}_t$}.}

In this paper we put into evidence a remarkable semigroup of 
homomorphisms of $( \cM_{+}, \boxtimes )$, where $\cM_{+}$ is the 
space of probability measures on $[ 0, \infty )$, and $\boxtimes$ is 
{\em free multiplicative convolution} -- the operation with 
measures from ${\mathcal M}_{+}$ which corresponds to the 
multiplication of free positive random variables (for a basic 
introduction to the ideas of free harmonic analysis, see for instance 
Chapter 3 of \cite{VDN}). 

For $\mu \in {\mathcal M}_{+}$ and $n \in {\mathbb N}$, the $n$-fold 
convolution $\mu \boxtimes \mu \boxtimes \cdots \boxtimes \mu$ is 
denoted by $\mu^{\boxtimes n}$. It turns out (\cite{BBercoviciMult}, Section 2) 
that for every $\mu \in {\mathcal M}_{+}$, the measures 
$\mu^{\boxtimes n}$ can be incorporated in a family 
$\{\mu^{\boxtimes t} \mid t \in [ 1, \infty ) \}$ such that
$\bigl( \, \mu^{\boxtimes s} \, \bigr) \boxtimes 
\bigl( \, \mu^{\boxtimes t} \, \bigr) \ = \  
\mu^{\boxtimes (s+t)}$, for all $s,t \geq 1$. In particular,
the extended family of $\boxtimes$-powers provides us with a 
continuous semigroup of homomorphisms for $\boxtimes$, consisting of 
the maps $\cM_{+} \ni \mu \mapsto \mu^{\boxtimes t} \in \cM_{+}$,
$t \geq 1$.

In this paper we show that, quite surprisingly, there is a simple 
formula which defines another semigroup of homomorphisms for 
$\boxtimes$, by using powers of two {\em additive} convolutions on 
${\mathcal M}_{+}$: the free additive convolution $\boxplus$,
and the Boolean convolution $\uplus$. (A brief review of $\boxplus$ 
and of $\uplus$ is made in Section 2 below.) This other semigroup of 
$\boxtimes$-homomorphisms consists of the maps
\begin{equation} \label{eqn:1.1}
{\mathcal M}_{+} \ni \mu \mapsto 
\Bigl( \, \mu^{\boxplus (1+t)} \, \Bigr)^{ \uplus (1/(1+t))}
\in {\mathcal M}_{+},
\ \  t \geq 0 .
\end{equation}

The formula shown in (\ref{eqn:1.1}) actually makes sense when 
$\mu$ belongs to the larger space ${\mathcal M}$ of all probability 
measures on ${\mathbb R}$ (without requiring that $\mu$ is supported 
on $[0, \infty )$). When moving to ${\mathcal M}$, one gets the issue 
that the convolution $\mu \boxtimes \nu$ isn't generally defined for 
arbitrary $\mu , \nu \in {\mathcal M}$. However, it is still possible 
to define $\mu \boxtimes \nu$ when $\mu \in {\mathcal M}$ and 
$\nu \in {\mathcal M}_{+}$, and we have the following theorem, which 
contains in particular the fact that the maps in (\ref{eqn:1.1}) 
form a semigroup of $\boxtimes$-homomorphisms of ${\mathcal M}_{+}$.

$\ $

{\bf Theorem 1.} 
{\em For every $t \geq 0$ one can define a one-to-one map 
${\mathbb B}_t : {\mathcal M} \to {\mathcal M}$ by
\begin{equation}  \label{eqn:1.4}
{\mathbb B}_t ( \mu ) = \Bigl( \, \mu^{\boxplus (1+t)} \, 
\Bigr)^{\uplus (1/(1+t))}, \ \ \mu \in {\mathcal M} .
\end{equation}
Every ${\mathbb B}_t$ is continuous with respect to the weak topology on
$\cM$, and carries $\cM_{+}$ into itself. 
Moreover, the maps $\{ \bB_{t} \mid t \geq 0 \}$ satisfy:
\begin{equation}  \label{eqn:1.6}
{\mathbb B}_s \circ {\mathbb B}_t =
{\mathbb B}_{s+t}, \ \  \forall \, s,t \geq 0,
\end{equation}
and for every $t \geq 0$ we have that 
\begin{equation}  \label{eqn:1.5}
{\mathbb B}_t ( \mu \boxtimes \nu ) =
{\mathbb B}_t ( \mu ) \boxtimes {\mathbb B}_t ( \nu), \ \ 
\forall \, \mu \in {\mathcal M}, \ \nu \in {\mathcal M}_{+} .
\end{equation}   }

$\ $

The semigroup of transformations $\{ {\mathbb B}_t \mid t \geq 0 \}$
has two other interesting features, which are discussed in the next
subsections 1.2 and 1.3.

\newpage

{\bf 1.2 Relation to \boldmath{$\boxplus$}-infinite divisibility.}

A probability measure $\mu \in {\mathcal M}$ is said to be 
{\em infinitely divisible} with respect to $\boxplus$ if for every 
$n \geq 1$ there exists $\mu_n \in {\mathcal M}$ such that 
$\mu_n^{\boxplus n} = \mu$; we will denote the set of probability 
measures which have this property by ${\mathcal M}_{\textrm{inf-div}}$. 
In terms of convolution powers with respect to $\boxplus$, the fact
that $\mu \in \cM_{\textrm{inf-div}}$ amounts to saying that 
$\mu^{\boxplus t}$ can be defined for every $t>0$ (in contrast with 
the situation of an arbitrary probability measure $\mu \in \cM$, for 
which the powers $\mu^{\boxplus t}$ can in general be defined only 
for $t \geq 1$). Infinite divisibility in free sense has a 
well-developped theory -- see 
section 2.11 of the survey \cite{V00}. An aspect of this theory which is 
of particular relevance for this paper is a special bijection
\begin{equation}  \label{eqn:1.7}
{\mathbb B} : {\mathcal M} \to {\mathcal M}_{\textrm{inf-div}},
\end{equation}
found by Bercovici and Pata (\cite{BPAnnals}, Section 6), in connection to 
their parallel study of infinite divisibility with respect to 
$\boxplus$ and to $\uplus$. We will refer to ${\mathbb B}$ as the 
{\em Boolean Bercovici-Pata bijection}. The transformations 
${\mathbb B}_t$ from our Theorem 1 connect to this as follows.

$\ $

{\bf Theorem 2.} {\em We have
\begin{equation}  \label{eqn:1.8} 
{\mathbb B}_1 ( \mu ) = {\mathbb B} ( \mu ), \ \ 
\forall \, \mu \in {\mathcal M},
\end{equation}
where the transformation 
${\mathbb B}_1 : {\mathcal M} \to {\mathcal M}$ is defined as in
Theorem 1 (by setting $t=1$ there), and ${\mathbb B}$ is 
the Boolean Bercovici-Pata bijection.  }

$\ $

A consequence of Theorem 2 and of Equation (\ref{eqn:1.5}) from 
Theorem 1 is that the Boolean Bercovici-Pata bijection is 
multiplicative with respect to $\boxtimes$. This phenomenon was 
observed, via combinatorial methods, to also hold for a 
multi-variable generalization of the Boolean Bercovici-Pata bijection 
which was recently studied in \cite{etac}.

Another consequence of Theorem 2, in conjunction with the semigroup 
property (\ref{eqn:1.6}) from Theorem 1 is that for $t \geq 1$ 
we have $\bB_t ( \cM ) = \bB ( \, \bB_{t-1} ( \cM ) \, ) \subseteq
{\mathcal M}_{\textrm{inf-div}}$. That is, we have the following 
corollary.

$\ $

{\bf Corollary.} {\em 
The probability measure ${\mathbb B}_t ( \mu )$ is infinitely 
divisible with respect to $\boxplus$, for every $t \geq 1$
and every $\mu \in {\mathcal M}$. }

$\ $

The statement of the above corollary can be sharpened by 
introducing a numerical quantity, defined as follows.

$\ $

{\bf Definition.} For $\mu \in \cM$ we denote
\begin{equation}  \label{eqn:1.71}
\phi ( \mu ) := \sup \{ t \in [ 0, \infty ) \mid 
\mu \in \bB_t ( \cM ) \}  \in [ 0, \infty ].
\end{equation}
We will call $\phi ( \mu )$ the {\em $\boxplus$-divisibility 
indicator} of $\mu$.

$\ $ 

In terms of the $\boxplus$-divisibility indicator, the statement of 
the preceding corollary gets translated into the fact 
that for $\mu \in \cM$ we have the equivalence 
\begin{equation}  \label{eqn:1.72}
\mu \in \cM_{\textrm{inf-div}} \ \Leftrightarrow \
\phi ( \mu ) \geq 1
\end{equation}
(see Proposition \ref{prop:5.4} below).
Thus if $\mu$ has $\phi ( \mu ) \geq 1$, then we are sure we can 
consider the $\boxplus$-power $\mu^{\boxplus t}$ for every $t>0$.
On the other hand, if $\mu \in \cM$ has $0 < \phi ( \mu ) < 1$ 
(and thus does not belong to $\cM_{\mathrm{inf-div}}$), then we can 
still take some subunitary $\boxplus$-powers of $\mu$ -- namely 
we can talk about $\mu^{\boxplus t}$ for every $t \ge 1- \phi (\mu )$;
see Remark \ref{rem:5.5}.3 below.

The values of $\phi ( \mu )$ for a few distributions of importance 
in free probability are listed in the next table.

\vspace{10pt}

\begin{tabular}{|c|c|||c|c|}  \hline
                     &          &             &        \\  \hline
Distribution $\mu$   &  $\phi ( \mu )$        &       
           Distribution $\mu$   &  $\phi ( \mu )$      \\        
                     &          &             &        \\  \hline
                     &          &             &        \\  \hline
                     &          &             &        \\
Symmetric Bernoulli distribution,   &      $0$    &    
Marchenko-Pastur distribution       &      $1$         \\
$\mu = \frac{1}{2} ( \delta_{-1} + \delta_1 )$  &   &    
of parameter 1,                                        \\
     &    &   $d \mu (x) = \frac{1}{2 \pi} \sqrt{(4-x)/x} \, dx$ on $[0,4]$  &  \\
                     &              &         &        \\  \hline
                     &              &         &        \\
Arcsine law of variance 1,          &                    $1/2$    &   
Cauchy distribution,                &      $\infty$    \\
$d \mu (x) = \frac{1}{\pi \sqrt{2-x^2}} \, dx$  
          on $[- \sqrt{2}, \sqrt{2}]$          &  &    
$d \mu (x) = \frac{1}{\pi (x^2 +1)} \, dx$        &    \\
                     &                        &   &    \\  \hline
                     &                        &   &    \\  
Standard semicircular distribution, &      $1$  & &    \\
$d \mu (x) = \frac{1}{2 \pi} \sqrt{4 - x^2} \, dx 
             \mbox{ on }  [-2,2]$          &    & &    \\ 
                     &                     &      &    \\  \hline
\end{tabular}

\vspace{10pt}

\begin{center}
{\bf Table 1.} 
$\phi ( \mu )$ for a few important distributions $\mu$.
\end{center}

$\ $

{\bf 1.3 Relation to free Brownian motion, and to the complex Burgers 
equation.}

For a probability measure $\mu \in {\mathcal M}$, we will use the 
notation $F_{\mu}$ for the reciprocal Cauchy transform of $\mu$. 
That is, $F_{\mu}$ is the analytic self-map of the upper half-plane 
$\bC^{+}$ defined as 
\[
F_{\mu} (z) := 1/G_{\mu} (z), \ \ \forall \, z \in {\mathbb C}^{+},
\]
where $G_{\mu}$ is the Cauchy transform of $\mu$ (a brief review of 
$F_{\mu}$ and of some of its basic properties appears in Section 2.2
below).

$\ $

{\bf Theorem 3.} {\em Let $\mu$ be in ${\mathcal M}$, and consider 
the function $h: ( 0, \infty ) \times \bC^{+} \to \bC$ defined by
\begin{equation}  \label{eqn:1.9}
h(t,z) = F_{\bB_t (\mu) } (z)-z, \ \  
\forall \, t > 0, \ \ \forall \, z\in {\mathbb C}^{+}.
\end{equation}
Then $h$ satisfies the complex Burgers equation,
\begin{equation}  \label{eqn:1.10}
\frac{\partial h}{\partial t}(t,z)=h(t,z)
\frac{\partial h}{\partial z}(t,z),\quad t>0,z\in \mathbb C^+.
\end{equation}  }

$\ $

The complex Burgers equation has previously appeared in free harmonic 
analysis in the work of Voiculescu \cite{V3}, in connection to the 
free Brownian motion started at a probability measure 
$\nu \in {\mathcal M}$ -- that is, in connection to the family of 
measures $\{ \nu \boxplus \gamma_t \mid t >0 \}$, where $\gamma_t$ is 
the semicircular distribution of variance $t$. More precisely, when 
one considers the Cauchy transforms of the measures in this family, 
it turns out that the function 
\begin{equation}  \label{eqn:1.11}
( 0, \infty ) \times {\mathbb C}^{+} \ni (t,z) \mapsto 
-G_{\nu \boxplus \gamma_t} (z) \in {\mathbb C}
\end{equation}
satisfies the complex Burgers equation (in exactly the form stated 
above in Equation (\ref{eqn:1.10})). These two occurrences of the 
complex Burgers equation (for the functions in (\ref{eqn:1.9}) and 
(\ref{eqn:1.11})) are in fact connected to each other, in the way 
described as follows. 

$\ $

{\bf Theorem 4.}
{\em Let $\nu$ be a probability measure 
in ${\mathcal M}$, and consider the analytic function $- G_{\nu}$
from ${\mathbb C}^{+}$ to itself. Then there exists a unique 
probability measure $\mu \in {\mathcal M}$ such that 
\begin{equation}  \label{eqn:1.12}
- G_{\nu} (z) = F_{\mu} (z) -z, \ \ z \in {\mathbb C}^{+}.
\end{equation}
Moreover, the relation (\ref{eqn:1.12}) between $\mu$ and $\nu$ is not
affected when $\nu$ evolves via the free Brownian motion, while $\mu$
evolves under the action of the transformations ${\mathbb B}_t$;
that is, we have that
\begin{equation}  \label{eqn:1.13}
- G_{\nu \boxplus \gamma_t} (z) = 
F_{{\mathbb B}_t ( \mu )} (z) -z, \ \ 
\forall \, t >0,  \ \ \forall \, z \in {\mathbb C}^{+}.
\end{equation}   }

$\ $

Clearly, Theorem 3 would follow from Theorem 4 and the 
corresponding result for the free Brownian motion, if it were true
that every analytic function from the collection 
$\{ F_{\mu} (z)-z \mid \mu \in {\mathcal M} \}$ can be put in the 
form $-G_{\nu}$ for some $\nu \in {\mathcal M}$; but this is not 
the case -- see Remark \ref{rem:4.3} below. Nevertheless,
Theorem 4 gives an interesting interpretation of the 
transformations ${\mathbb B}_t$, as a way of 
re-casting the free Brownian motion where the resulting 
process becomes multiplicative with respect to $\boxtimes$.

$\ $

{\bf 1.4 Further remarks, and organization of the paper.}

In \cite{IJM}, using combinatorial methods, we have found
multi-variable analogues to several results we present in
this paper. Most notably, \cite{IJM} provides a multi-variable
analogue of Theorem 4, and an operatorial model for the correspondence
described in this theorem. However, the methods of \cite{IJM} have
the important draw-back that they only apply to distributions
with compact support. All the results of the present paper use
complex analytic tools, and apply to arbitrary Borel probability measures on the real line.

Let us now give a brief
outline of the paper's organization:
After a review of background material in Section 2, we give the 
proofs of Theorems 1 and 2 in Section 3. What stands behind Theorem 1 
are two interesting connections which the Boolean convolution 
$\uplus$ turns out to have with the operations $\boxplus$ and 
$\boxtimes$ from free probability.
The first connection is a kind of {\em commutation relation} between  
the convolution powers with respect to $\boxplus$ and to $\uplus$: a 
probability measure  of the form $( \mu^{\boxplus p} )^{\uplus q}$ can 
also be written as $( \mu^{\uplus q'} )^{\boxplus p'}$ for some new 
exponents $p', q'$, where $p', q'$ are given explicitly in terms 
of $p$ and $q$ -- see Proposition \ref{prop:3.1} below. The second 
connection is a kind of ``distributivity'' which involves 
$\boxtimes$ and a fixed convolution power with respect to $\uplus$:
\begin{equation}   \label{eqn:1.130}
( \mu^{\uplus t} ) \boxtimes ( \nu^{\uplus t} ) = 
( \mu \boxtimes \nu )^{\uplus t} \circ D_{1/t}, 
\ \  \forall \, t > 0, \ \forall \, \mu , \nu \in \cM .
\end{equation}
In (\ref{eqn:1.130}), ``$D_{1/t}$'' stands for the natural operation 
of dilation (by a factor of $1/t$) for probability measures on $\bR$. 
Together with an analogous distributivity formula involving $\boxtimes$ 
and a convolution power of $\boxplus$, this relation explains why 
each of the transformations $\bB_t$ from Theorem 1 is a homomorphism 
with respect to $\boxtimes$. See Proposition \ref{prop:3.4} and 
Remark \ref{rem:3.5} below.

In Section 3 we also point out that the transformations $\bB_t$ 
can be very nicely described by using Voiculescu's $S$-transform 
(see Remark \ref{rem:3.6}). Theorem 2 can be easily proved as an 
application of this (see Remark \ref{rem:3.7}).

Section 4 is devoted to the relation with the free Brownian motion, 
and to the proofs of Theorems 3 and 4. A very nice example of 
process $\{ \bB_t ( \mu ) \mid t \geq 0 \}$ which can be described 
explicitly by using Theorem 4 is the one started at the symmetric 
Bernoulli distribution; this process turns out to go both through 
the arcsine law and through the standard semicircle law (see Example 
\ref{ex:4.5}).

Finally, in Section 5 we discuss a few miscellaneous facts related to the 
transformations $\id_t$. In Proposition \ref{prop:5.1} and Corollary 
\ref{cor:5.2} we describe atoms and regularity for the measures 
$\id_t(\mu)$ and the convergence
of $\id_t(\mu)$ when $t$ tends to zero. In Remarks \ref{rem:5.5} and 
\ref{rem:5.6} we prove some basic properties of the 
$\boxplus$--divisibility indicator $\phi ( \mu )$ introduced 
in Section 1.2, including the verification for the values of 
$\phi ( \mu )$ that were listed in Table 1.

\section{Background and notations}

\setcounter{equation}{0}

$\ $

{\bf 2.1 The convolution operations 
${\mathbf \boxplus ,  \boxtimes ,  \uplus .}$ }

Same as in the introduction, we use the notation $\cM$ for the set 
of Borel probability measures on $\bR$, and the notation $\cM_+$ 
for the set of all probability measures $\mu \in \cM$ with the 
property that $\mu([0,+\infty))=1.$

In the literature on non-commutative probability one encounters
several ``convolution'' operations for probability distributions 
in $\cM$, which are defined to reflect operations with 
non-commutative random variables.

The operations $\boxplus$ and $\boxtimes$ are from free probability
theory. They are defined in order to reflect the addition 
and respectively the multiplication of free random variables; we 
refer to Chapter 3 of \cite{VDN} or Chapters 2, 3 in \cite{V00} for 
a precise description of how this goes. In this paper we will not pursue 
the approach to $\boxplus$ and $\boxtimes$ in terms of free random 
variables, but we will rather use the analytic function theory 
developped in order to deal (and do computations) 
with these operations. In particular, in the next subsection 2.2 
we will review the ``transforms'' that are mainly used to study 
these operations: the $R$-transform for $\boxplus$,
and the $S$-transform for $\boxtimes$.

There is a third convolution operation which appears in this paper,
the Boolean additive convolution $\uplus$. This comes from the world
of random variables that are ``Boolean independent'', and reflects 
the addition of such variables. One of the points emphasized by 
this paper is that $\uplus$ has nice relations with $\boxplus$ and
$\boxtimes$, and, surprisingly, has a role to play in free 
probability, in connection to infinite divisibility with respect 
to $\boxplus$. In the next subsection we will also review how 
$\uplus$ is handled by using complex function theory -- the 
analytic function theory for $\uplus$ is in fact simpler than the
one required for $\boxplus$ and $\boxtimes$.

$\ $

{\bf 2.2 A glossary of transforms.}

Here we collect some basic formulas concerning the various kinds of
transforms used in non-commutative probability, and which appear in 
this paper. We are not aiming to a self-contained
presentation of the transforms, for the most part we will only 
state the formulas and properties that we need, and indicate 
references for them.

\vspace{6pt}

$1^o$ {\bf Cauchy transform and reciprocal Cauchy transform.}

The {\bf Cauchy transform} of a probability measure $\mu \in \cM$ 
is the analytic function $G_{\mu}$ defined by
\begin{equation}  \label{eqn:2.1}
G_{\mu} (z)=\int_\mathbb R\frac{d \mu (s)}{z-s},
\quad z\in\mathbb C \setminus \bR.
\end{equation}
The {\bf reciprocal Cauchy transform} $F_{\mu}$ is defined by
\begin{equation}  \label{eqn:2.2}
F_{\mu} (z)= 1/ G_{\mu} (z), \ \ z \in \bC \setminus \bR.
\end{equation}
It can be easily checked that $G_{\mu}$ maps $\bC^+$ to $\bC^-$;
as a consequence of this, $F_{\mu}$ can be viewed 
as an analytic self-map of $\bC^+$.

We will denote by $\mathcal F$ the set of analytic self-maps of 
$\bC^+$ which can appear as $F_{\mu}$ for some $\mu \in \cM$. 
One has a very nice intrinsic description for the functions in 
$\mathcal F$, namely 
\begin{equation}  \label{eqn:2.3}
\mathcal F = \Bigl\{ F : \bC^+ \to \bC^+ \mid \mbox{ $F$ is 
analytic and } 
\lim_{y \to \infty} \frac{F(iy)}{iy} = 1 \Bigr\} .
\end{equation}
For a function $F \in \mathcal F$, the limit $\lim F(z)/z = 1$ 
holds in fact under the weaker condition that $z$ converges 
non-tangentially to $\infty$ (i.e. $z \to \infty$ in 
an angular domain of the form 
$\{ z \in \bC^+ \mid |\Re (z)| < c \cdot \Im (z) \}$, 
for some $c>0$).

Another fact worth recording is that a function $F \in {\mathcal F}$ 
always increases the imaginary part: 
$\Im F_\rho(z) \geq \Im z$ for all $z \in \bC^+$.
Moreover, if there exists $z_0\in\mathbb C^+$ such that 
$\Im F_\rho(z_0) = \Im z_0,$ then the equality holds for all 
$z\in\mathbb C^+$ and $\rho$ is a point mass.
Proofs of all these facts can be found in \cite{Achieser}. For a 
nice review, one can also consult Section 2 in Maassen's paper
\cite{M92} or Section 5 of \cite{BVIUMJ}.

\vspace{6pt}

$2^o$ {\bf $R$-transform.}

Let $\mu$ be a probability measure in $\cM$, and consider its 
Cauchy transform $G_{\mu}$. It can be proved
(\cite{Voiculescu1},\cite{VDN},\cite{BVIUMJ}) that 
the composition inverse $G_\mu^{-1}$ of $G_\mu$ is defined on 
a truncated angular domain of the form 
\[
\Bigl\{ z \in \bC^+ \mid | \Re (z) | < c \cdot \Im (z), \ 
|z|> M \Bigr\} , 
\]
for some $c, M>0$. For $z$ in this domain we define 
\[
\mathcal R_\mu(z)=G_\mu^{-1}(z)-\frac1z, \quad{\rm and}
\quad R_\mu(z)=z\mathcal R_\mu(z).
\]
The function $R_{\mu}$ is called the {\bf $R$-transform} of $\mu$.
(Note that in the free probability literature the function 
$\mathcal R_{\mu}$ also goes under the same name, of $R$-transform
of $\mu$; in this  paper we will only work with $R_{\mu}$.)

The $R$-transform is the linearizing transform for the operation
$\boxplus$. That is, for every $\mu , \nu \in \cM$ we have
\begin{equation} \label{eqn:2.4}
R_{\mu \boxplus \nu } (z) = R_{\mu} (z) + R_{\nu} (z),
\end{equation}
with $z$ running in truncated angular domain where all of 
$R_{\mu} (z)$, $R_{\nu} (z)$ and $R_{ \mu \boxplus \nu } (z)$ 
are defined.

\vspace{6pt}

$3^o$ {\bf $S$-transform.}

Let $\mu \in \cM$ be a probability measure in $\cM$ which has 
compact support and has the first moment $\int t \ d \mu (t)$
different from $0$. Consider the moment generating series of $\mu$,
\[
\psi_{\mu} (z) = \sum_{n=1}^{\infty} m_n z^n,
\]
defined on a neighborhood of $0$, and where for every $n \geq 1$
we denote $m_n = \int t^n \ d \mu (t)$ (the moment of order $n$ 
of $\mu$). Then $\psi_{\mu}$ is invertible under composition on 
a sufficiently small disc centered at $0$, and it makes sense 
to define
\begin{equation} \label{eqn:2.5}
S_\mu(z)=\frac{z+1}{z}\psi_\mu^{-1}(z), \ \ 
\mbox{ for $|z|$ sufficiently small.}
\end{equation}
$S_{\mu}$ is called the {\bf $S$-transform} of $\mu$.

An equivalent way of defining the $S$-transform is by relating 
it directly to the $R$-transform $R_{\mu}$. Indeed, for $\mu$ as 
above it is easily seen that $R_{\mu}$ is defined and invertible 
under composition on a small disc centered at $0$, and it turns 
out that we have 
\begin{equation}  \label{eqn:2.6}
S_\mu(z)=\frac{1}{z}R_\mu^{-1}(z),
\end{equation}
again holding for $z \in \bC$ such that $|z|$ is sufficiently
small.

It can be proved (\cite{V2}) that the $S$-transform is multiplicative 
in for $\boxtimes$, in the sense that if $\mu \in \cM$ and 
$\nu \in \cM_+$ have compact support and the first moment different 
from zero, then 
\begin{equation}  \label{eqn:2.7}
S_{\mu \boxtimes \nu}(z) = S_{\mu}(z)  S_{\nu}(z).
\end{equation}

\vspace{6pt}

$4^o$ {\bf $\eta$-transform.}

For a probability measure $\mu \in \cM$ we denote
\begin{equation}  \label{eqn:2.8}
\psi_\mu (z)=\int_\mathbb R\frac{sz}{1-sz}d\mu(s),
\quad z \in \mathbb C \setminus \bR .
\end{equation}
We will refer to $\psi_{\mu} (z)$ as the {\em moment generating 
function} of $\mu$. (In the case when $\mu$ has compact support, 
it is easily seen that the integral formula given in 
(\ref{eqn:2.8}) matches the series expansion, also denoted by 
$\psi_{\mu}$, which appeared in the above discussion about the 
$S$-transform.)

The {\bf $\eta$-transform} of $\mu$ is then defined in terms of the 
moment generating function $\psi_{\mu}$ by the formula
\begin{equation}   \label{eqn:2.9}
\eta_\mu (z)=\frac{\psi_\mu(z)}{ 1+\psi_\mu(z) },
\ \ z \in \bC \setminus \bR .
\end{equation}
Note that the denominator of the fraction on the right-hand side 
of Equation (\ref{eqn:2.9}) is always different from zero, due to 
the fact (immediately seen from the definitions) that we have
\begin{equation}  \label{eqn:2.10}
1 + \psi_{\mu} (z) = \frac{1}{z} G_{\mu} (1/z) \neq 0, 
\ \ \forall \, z \in \bC \setminus \bR .
\end{equation}
Let us also record here that $\eta_{\mu}$ can be directly obtained 
from the reciprocal Cauchy transform $F_{\mu}$ via the formula
\begin{equation}  \label{eqn:2.11}
\eta_{\mu} (z) = 1 - z F_{\mu} (1/z), \ \ z \in \bC \setminus \bR 
\end{equation}
(which is obtained by combining the preceding two Equations 
(\ref{eqn:2.9}) and (\ref{eqn:2.10})).

The $\eta$-transform plays a role in the study of the Boolean 
convolution $\uplus$. Indeed, it is known \cite{SpeicherWoroudi} that $\uplus$ 
can be neatly described by using reciprocal Cauchy transforms: for 
$\mu, \nu \in \cM$ we have
\begin{equation}  \label{eqn:2.110}
F_{\mu \uplus \nu} (z) = F_{\mu} (z) + F_{\nu} (z) - z, \ \ 
\forall \, z \in \bC \setminus \bR.
\end{equation}
This amounts to saying that ``$F_{\mu} (z) -z$ is a linearizing 
transform for $\uplus$'' (in the sense that the function 
``$F(z)-z$'' calculated for $\mu \uplus \nu$ is the sum of the 
corresponding functions calculated for $\mu$ and $\nu$). But now,
thanks to the direct connection between $\eta$-transform and 
reciprocal Cauchy transform which was recorded in (\ref{eqn:2.11}), 
one immediately sees that the latter linearization formula 
(\ref{eqn:2.110}) is equivalent to
\begin{equation}  \label{eqn:2.120}
\eta_{\mu \uplus \nu} (z) = \eta_{\mu} (z) + \eta_{\nu} (z), 
\ \ \forall \, \mu , \nu \in \cM, \ \forall \, z \in \bC^+.
\end{equation}
Thus we can conclude that the $\eta$-transform also is ``a 
linearizing transform for $\uplus$''.

\vspace{6pt}

$5^o$ {\bf $\Sigma$-transform.}

By comparing the above Equations (\ref{eqn:2.4}) and 
(\ref{eqn:2.120}), one could say that the $\eta$-transform is an 
analogue of the $R$-transform, living in the parallel world of 
Boolean probability. But the $\eta$-transform also has a direct 
connection (which is more than a ``Boolean vs. free'' analogy) with 
the $R$-transform. This connection is in terms of the composition 
inverses for $R$ and $\eta$, and is best put into evidence by 
introducing the $\Sigma$-transform of $\mu$.

Let us assume that $\mu \in \cM$ has compact support and has the 
first moment different from $0$. It is easily seen that in this
case the $\eta$-transform $\eta_{\mu}$ is defined and invertible 
under composition on a small disc centered at $0$, and it thus
makes sense to define the {\bf $\Sigma$-transform} 
$\Sigma_{\mu}$ of $\mu$ by the formula
\begin{equation}  \label{eqn:2.130}
\Sigma_\mu(z) =\frac{1}{z}\eta_\mu^{- 1}(z),
\mbox{ for $|z|$ sufficiently small.}
\end{equation}
Note that this is very similar to the formula used when one 
defines the $S$-transform in terms of the $R$-transform
(see Equation (\ref{eqn:2.6}) above). But more than having an analogy
between how $S_{\mu}$ and $\Sigma_{\mu}$ are defined, it turns 
out that we have a direct relation between them, namely:
\begin{equation}  \label{eqn:2.140}
\Sigma_\mu(z)=S_\mu\left(\frac{z}{1-z}\right),
\mbox{ $|z|$ sufficiently small.}
\end{equation}
Equation (\ref{eqn:2.140}) is obtained by taking inverses under 
composition in Equation (\ref{eqn:2.9}) (this gives us that 
$\eta_{\mu}^{-1} (z) = \psi_{\mu}^{-1} (z/(1-z))$, for $|z|$ sufficiently 
small), and then by invoking the formulas (\ref{eqn:2.5}) and 
(\ref{eqn:2.130}) which were used to define $S_{\mu}$ and 
respectively $\Sigma_{\mu}$.

${ }$From Equation (\ref{eqn:2.140}) it is immediate that the 
$\Sigma$-transform has a multiplicativity property with respect to
$\boxtimes$, analogous to the one enjoyed by the $S$-transform:
\begin{equation}  \label{eqn:2.150}
\Sigma_{\mu\boxtimes\nu}(z)=\Sigma_\mu(z)\Sigma_\nu(z),
\end{equation}
whenever $\mu \in \cM$ and $\nu \in \cM_+$ have compact support
and have first moment different from zero.

$\ $

{\bf 2.3 Remark.} One should keep in mind that a probability measure
$\mu \in \cM$ is uniquely determined by any of the transforms 
$(G_{\mu}, F_{\mu}, R_{\mu}, \ldots \, )$ reviewed in the above 
glossary, and which are defined for $\mu$. So, for instance, if 
$\mu , \nu \in \cM$ are such that $F_{\mu} = F_{\nu}$, then it
follows that $\mu = \nu$. Or: if $\mu, \nu \in \cM$ have compact
support and first moment different from $0$, and if we know that 
$S_{\mu}$ and $S_{\nu}$ coincide on a neighborhood of $0$, then 
we can infer that $\mu = \nu$. The reason for this is that each 
of the transforms considered in the above glossary determines the 
Cauchy transform of the measure; and a measure $\mu \in \cM$ can
always be retrieved from its Cauchy transform $G_{\mu}$, by a
procedure called ``the Stieltjes inversion formula'' (see e.g. 
\cite{Achieser}).

$\ $

In the remaining part of this section we review some facts 
about convolution powers and infinite divisibility that are 
used later on in the paper.

$\ $

{\bf 2.4 Convolution powers with respect to $\boxplus$, $\uplus$.}

$1^o$ For $\mu \in \cM$ and a positive integer $n$, one denotes the 
$n$-fold convolution $\mu \boxplus \cdots \boxplus \mu$ by 
$\mu^{\boxplus n}$. The probability measure $\mu^{\boxplus n}$ is 
very nicely characterized in terms of $R$-transforms, via the
formula $R_{\mu^{\boxplus} n} (z) = n \cdot R_{\mu} (z)$.
It turns out that the latter formula can be extended to the case
when $n$ is not an integer. More precisely, for every $\mu \in \cM$ 
and $t \in [ 1, \infty )$ there exists a probability measure 
$\mu_t \in \cM$ so that 
\begin{equation} \label{eqn:2.160}
R_{\mu_t}(z)= tR_\mu(z),
\end{equation}
with $z$ running in a truncated angular domain where both sides of 
the equation are defined. The existence of $\mu_t$ was first 
observed in \cite{NS} in the case when $\mu$ has compact support, 
and then extended to arbitrary $\mu \in \cM$ in \cite{BB}.

The measure $\mu_t$ appearing in (\ref{eqn:2.160}) (which is
uniquely determined by the prescription for the $R$-transform 
$R_{\mu_t}$) is called the $t$-th {\bf convolution power} of $\mu$ 
with respect to $\boxplus$, and denoted by $\mu^{\boxplus t}$. 
It is immediate that we have
\[
\mu^{\boxplus t} \boxplus \mu^{\boxplus s} = 
\mu^{\boxplus t+s}, \  \ \forall \, s,t \in [1,+\infty).
\]

In the following sections we shall use some other properties of the 
$\boxplus$-convolution powers.
It is an immediate consequence of the operatorial realization of $\mu^{\boxplus t}$ given 
in \cite{NS} that if $\mu\in\mathcal M$ has compact support, then so does $\mu^{\boxplus t}$
for any $t\ge1$, and if the support of $\mu$ is included in $[0,+\infty),$ then so is the support 
of $\mu^{\boxplus t}$ for any $t\ge1$. The operation of taking convolution powers with respect to $\boxplus$ is
well behaved with respect to the usual topologies on $\mathcal M$ and the real line:
it follows from equation (\ref{eqn:2.160}) and the characterization of continuity in terms of the 
$R$-transform given in \cite{BVIUMJ} that the correspondence $\mu\mapsto\mu^{\boxplus t}$
is continuous in the weak topology for any $t,$ and the correspondence $[1,+\infty)\ni 
t\mapsto\mu^{\boxplus t}\in\mathcal M$
is continuous for every $\mu\in\mathcal M$, 
where $[1,+\infty)$ is considered with the usual topology and $\mathcal M$ is endowed with the
weak topology. As a consequence of the two properties above, we 
easily see that if $\mu \in \cM_+$ then $\mu^{\boxplus t} \in \cM_+$, 
$\forall \, t \geq 1$.

\vspace{6pt}

$2^{o}$ We now do the same kind of discussion as above, but in 
connection to the operation of Boolean convolution $\uplus$.
For $\mu \in \cM$ and a positive integer $n$, one denotes the 
$n$-fold convolution $\mu \uplus \cdots \uplus \mu$ by 
$\mu^{\uplus n}$. From the discussion about the linearizing of
$\uplus$ (see Equation (\ref{eqn:2.110}) in the glossary of 
transforms) we see that $\mu^{\uplus n}$ is nicely characterized 
in terms of its reciprocal Cauchy transform, via the formula
\[
F_{\mu^{\uplus n}} (z) = n F_{\mu} (z) + (1-n) z, \ \ 
z \in \bC^+ .
\]
This formula can be extended to the case when the exponent $n$ is not an 
integer; in fact, it turns out that the convolution power $\mu^{\uplus t}$ 
can be defined for every $\mu \in \cM$ and every $t>0$ (we no longer have 
the restriction ``$t \geq 1$'' which we had in the above discussion 
about $\boxplus$). That is, for every $\mu \in \cM$ and every $t>0$, the 
$\uplus$-convolution power $\mu^{\uplus}$ is defined \cite{SpeicherWoroudi} as the 
unique probability measure in $\cM$ whose reciprocal Cauchy transform 
satisfies
\begin{equation}  \label{eqn:2.170}
F_{\mu^{\uplus t}} (z) = t F_{\mu} ( z) + (1-t)z, \ \ z \in \bC^+.
\end{equation}

It follows easily from properties of the functions $F_\mu$ that convolution powers with respect to
$\uplus$ enjoy the same properties as the ones 
mentioned for $\boxplus$ at the end of part $1^o$ of this remark.

$\ $

{\bf 2.5 Infinite divisibility with respect to $\boxplus$, $\uplus$.}

$1^o$ A probability measure $\mu \in \cM$ is said to be 
infinitely divisible with respect to $\boxplus$ if for every 
$n \geq 1$ there exists $\mu_n \in \cM$ such that 
$\mu_n^{\boxplus n} = \mu$. The set of 
probability measures $\mu \in \cM$ which are 
$\boxplus$-infinitely divisible will be denoted in this paper by 
$\cM_{\mathrm{inf-div}}$. It is easily seen that for a measure 
$\mu \in \cM_{\mathrm{inf-div}}$ one can define the 
$\boxplus$-convolution powers $\mu^{\boxplus t}$ in exactly the 
same way as described in Equation (\ref{eqn:2.160}), and where 
now $t$ can be an arbitrary number in $( 0, \infty )$. (And 
conversely, if $\mu \in \cM$ has the property that 
$\mu^{\boxplus t}$ is defined for all $t>0$, then it is immediate 
that $\mu \in \cM_{\mathrm{inf-div}}$.)

Infinite divisibility with respect to $\boxplus$ was first studied
in \cite{Voiculescu1}, then in \cite{BVIUMJ}. It was observed there 
that $\boxplus$-infinite divisibility is very nicely described in 
terms of the $R$-transform: given $\mu \in \cM$, one has that 
$\mu$ is $\boxplus$-infinitely divisible if and only if the 
$R$-transform $R_{\mu}$ can be extended analytically to all of
$\bC \setminus \bR$.

$2^o$ In the same vein as above, one could consider the parallel 
concept of infinite divisibility with respect to $\uplus$. But here 
the situation turns out to be much simpler. Indeed, as already 
pointed out in part $2^o$ of the preceding subsection, we have
that {\em every} $\mu \in \cM$ is infinitely divisible with respect 
to $\uplus$. (In particular, no new notation is needed for the 
set of measures in $\cM$ which are $\uplus$-infinitely divisible.)

$\ $

{\bf 2.6 The Boolean Bercovici-Pata bijection.}

In the paper \cite{BPAnnals}, Bercovici and Pata have proved the 
existence of a strong connection between free, Boolean, and 
classical infinite divisibility. We reproduce here the result 
from \cite{BPAnnals} which is relevant for the present paper. Let 
$( \mu_n )_{n=1}^{\infty}$ be a sequence in $\mathcal M$, and let
$k_1 < k_2 < \cdots < k_n < \cdots$ be 
a sequence of positive integers. Then (as proved in \cite{BPAnnals},
Theorem 6.3) the following statements are equivalent:
\begin{enumerate} 
\item[{\rm(1)}] The sequence 
$\underbrace{ \mu_n\boxplus\mu_n\boxplus\cdots \boxplus 
\mu_n }_{ k_n \ {\rm times} }$ converges weakly to 
a probability measure $\nu \in \cM$.

\item[{\rm(2)}] The sequence 
$\underbrace{ \mu_n \uplus \mu_n \uplus \cdots \uplus
\mu_n }_{ k_n \ {\rm times} }$ converges weakly to 
a probability measure $\mu \in \cM$.
\end{enumerate}

\noindent
Moreover, suppose that the statements (1) and (2) are both true. 
Then the limit $\nu$ from (1) is $\boxplus$-infinitely divisible, 
and we have the following relation between $\mu$ and $\nu$:
\begin{equation}  \label{eqn:2.180}
z-F_{\mu}(z)=zR_{\nu}(1/z),\quad z\in\mathbb C\setminus\mathbb R.
\end{equation}
Note that the right-hand side of the above equation makes indeed sense 
because, as mentioned in the preceding subsection, the $R$-transform 
of a $\boxplus$-infinitely divisible distribution extends 
analytically to all of $\bC \setminus \bR$.

And finally (this is also part of Theorem 6.3 in \cite{BPAnnals}),
the correspondence $\mu \mapsto \nu$ with $\mu, \nu$ as in Equation
(\ref{eqn:2.180}) is a bijection between $\cM$ and 
$\cM_{\mathrm{inf-div}}$. This correspondence will be called the
{\bf Boolean Bercovici-Pata bijection}, and will be
denoted by $\bB$. 

Let us record here that an alternative form of the Equation 
(\ref{eqn:2.180}) describing the Boolean Bercovici-Pata bijection is
\begin{equation}  \label{eqn:2.20}
R_{\bB ( \mu )} (z) = \eta_{\mu} (z), \ \ z \in \bC \setminus \bR .
\end{equation}
This is obtained by replacing $z$ with $1/z$ in 
(\ref{eqn:2.180}), and by invoking the formula (\ref{eqn:2.11})
which connects the $\eta$-transform to the reciprocal Cauchy transform.

Another useful reformulation of the Equation 
(\ref{eqn:2.180}) describing the bijection $\bB$ is obtained by 
using the $S$-transform and the $\Sigma$-transform. Let us suppose
that $\mu \in \cM$ has compact support and has first moment different 
from $0$. Then the equality from (\ref{eqn:2.20}) can be extended 
to a small disc centered at $0$, where we can also perform 
inversion under composition for the functions on its two sides. 
Due to the analogy between how $S$ and $\Sigma$ were defined in terms 
of $R$ and respectively $\eta$ (see Equations (\ref{eqn:2.6}) and 
(\ref{eqn:2.130}) in the above glossary of transforms), we thus 
arrive to the fact that
\begin{equation}  \label{eqn:2.21}
S_{\bB ( \mu )} (z) = \Sigma_{\mu} (z), \ \ 
\mbox{for $|z|$ sufficiently small.}
\end{equation}

We conclude this subsection by recording a few other properties of the 
bijection $\bB$, which were established in \cite{BPAnnals} or follow 
easily from the arguments presented there.
The bijection $\id$ turns out to be a weakly continuous isomorphism from $(\mathcal M, \uplus)$
onto $(\mathcal M_{\mathrm{inf-div}},\boxplus)$. 
Moreover, as an easy consequence of relation (\ref{eqn:2.20}),
if $\mu$ has compact support, then so does $\id(\mu).$

$\ $

{\bf 2.7 Subordination.}

Let $\mu$ and $\sigma$ be two probability measures in $\cM$. One says
that $F_\sigma$ is subordinated to $F_\mu$ (abberviated in this paper as ``$\sigma$ 
is {\bf subordinated} to $\mu$'') if there exists a 
function $\omega \in \mathcal F$ (with $\mathcal F$ as described in 
Equation (\ref{eqn:2.3}) above), such that 
\[
F_{\sigma} = F_{\mu} \circ \omega .
\]
If it exists, this function $\omega \in \mathcal F$ is uniquely 
determined, and is called the {\bf subordination function} of 
$\sigma$ with respect to $\mu$. An important phenomenon in the study
of the operation of free additive convolution $\boxplus$ is that 
for every $\mu, \nu \in \cM$, the convolution $\mu \boxplus \nu$ is 
subordinated with respect to $\mu$ and with respect to $\nu$.
This phenomenon has been first observed in \cite{V3} and then proved in 
full generality in \cite{Biane1}.

In the same vein, one has that a $\boxplus$-convolution power 
$\mu^{\boxplus t}$ is always subordinated with respect to $\mu$, 
for every $\mu \in \cM$ and every $t \geq 1$. This fact appears
in \cite{BB}. In the same paper it is also pointed out that the 
subordination function $\omega$ of $\mu^{\boxplus t}$ with respect 
to $\mu$ can be given by a ``direct'' formula (not involving 
composition) in terms of the reciprocal Cauchy transform of 
$\mu^{\boxplus t}$, namely
\begin{equation}  \label{eqn:2.22}
\omega (z) = \frac{1}{t} z + \Bigl( 1 - \frac{1}{t} \Bigr)
F_{\mu^{\boxplus t}} (z), \ \ z \in \bC^+ .
\end{equation}
We will repeatedly use this formula in what follows, rewritten in 
order to express $F_{\mu^{\boxplus t}}$ in terms of $\omega$, 
when it thus says that
\begin{equation}  \label{eqn:2.23}
F_{\mu^{\boxplus t}} (z) = \frac{t\omega(z)-z}{t-1}, \ \ z \in \bC^+ 
\end{equation}
(where $\mu$, $t$ and $\omega$ are the same as in (\ref{eqn:2.22})). 

On the other hand let us note that the above Equation
(\ref{eqn:2.22}) can be also put in the form
\begin{equation}   \label{eqn:2.235}
\omega (z) = \frac{1}{t} z + \Bigl( 1 - \frac{1}{t} \Bigr)
F_{\mu} ( \omega (z) ). 
\end{equation}
The latter formula can be viewed as a functional equation, for 
which it is known that $\omega$ is the only solution belonging to
the set of analytic maps $\mathcal F$ from Equation (\ref{eqn:2.3}).

Another benefit of Equation (\ref{eqn:2.235}) is that we can use 
it in order to write $z$ in terms of $\omega (z)$:
\begin{equation}  \label{eqn:2.24}
z = t \omega (z) + (1-t) F_{\mu} ( \omega (z) ), \ \ z \in \bC^+ ,
\end{equation}
and the latter means that the function 
\begin{equation}  \label{eqn:2.25}
H : \bC^+ \to \bC, \ \ \ H(w) = t w + (1-t) F_{\mu} ( w ), 
\end{equation}
is a left-inverse for $\omega$. It was in fact observed in 
\cite{BB} (and will be used in this paper too) that an 
equality of subordination functions of the kind discussed above 
is equivalent to the equality of their left-inverses defined as 
in (\ref{eqn:2.25}). More precisely: let $\mu, t, \omega, H$ be 
as above, and let $\widetilde{\mu}$, $\widetilde{t}$, 
$\widetilde{\omega}$, $\widetilde{H}$ be another set of data 
given in the same way; then we have that
\begin{equation}  \label{eqn:2.26}
\omega = \widetilde{\omega} \ \Leftrightarrow \
H = \widetilde{H}.
\end{equation}

Finally, let us record here one more formula concerning subordination 
functions, which will be invoked in Section 4 below. This formula 
is not about $\boxplus$-convolution powers, but rather concerns the 
free Brownian motion, i.e. the free additive convolution with a 
semicircular distribution. For $t>0$, let us denote by 
$\gamma_t$ the centered semicircular distribution of variance $t$,
$d \gamma_t (x) = \frac{1}{2 \pi t} \sqrt{4t -x^2} \, dx$ on 
$[-2 \sqrt{t}, 2 \sqrt{t} ]$. Let $\mu$ be a probability measure 
in $\cM$, let $t$ be in $(0, \infty )$, and let $\omega$ be the 
subordination function of $\mu \boxplus \gamma_t$ with respect to 
$\mu$. Then one has a ``direct'' formula (not involving compositions)
which relates $\omega$ to the Cauchy transform of 
$\mu \boxplus \gamma_t$, namely 
\begin{equation} \label{eqn:2.27} 
\omega (z) = z - t G_{\mu \boxplus \gamma_t} (z), \ \ z \in \bC^+ . 
\end{equation}
This formula was observed in \cite{Biane2} (it is obtained by 
putting together the statements of Lemma 4 and Proposition 2 of that 
paper).

\section{Proof of Theorems 1 and 2}

The fact that the 
\setcounter{equation}{0}
transformations $\{ \bB_t \mid t \geq 0 \}$ form a semigroup under 
composition will follow from a ``commutation relation'', stated in 
the next proposition, concerning the convolution powers with respect
to $\boxplus$ and to $\uplus$.

\begin{proposition}   \label{prop:3.1}
Let $p,q$ be two real numbers such that $p \geq 1$ and 
$q > (p-1)/p$. We have
\begin{equation} \label{eqn:3.1}
\Bigl( \, \mu^{\boxplus p} \, \Bigr)^{\uplus q} = 
\Bigl( \, \mu^{\uplus q'} \, \Bigr)^{\boxplus p'}, \ \ 
\forall \, \mu \in {\mathcal M},
\end{equation}
where the new exponents $p'$ and $q'$ are defined by
\begin{equation} \label{eqn:3.2}
p' := pq/(1-p+pq), \ \ q' := 1-p+pq
\end{equation}
(note that $p' \geq 1$ and $q' > 0$, thus the convolution powers 
appearing on the right-hand side of (\ref{eqn:3.1}) do indeed 
make sense).  
\end{proposition}

\begin{proof} 
If $p=1$, then we have $p'=1$ and $q'=q$, and both sides of 
(\ref{eqn:3.1}) are equal to $\mu^{\uplus q}$. Throughout the 
remaining of the proof we will assume that $p>1$ (which also implies 
that $p'>1$, thus allowing divisions by $p-1$ and $p'-1$ in our 
calculations).

We will prove the equality (\ref{eqn:3.1}) by showing
that the probability measures appearing on its two sides have 
the same reciprocal Cauchy transform. In order to do this we will
take advantage of the specific formulas we have for reciprocal 
Cauchy transforms, when we look at convolution powers for $\uplus$ and 
for $\boxplus$. For $\uplus$ we will simply use the formula 
(\ref{eqn:2.170}) which defined $\mu^{\uplus t}$ in the above 
subsection 2.4. The convolution powers with respect to $\boxplus$ 
are a bit more complicated: we will handle their reciprocal 
Cauchy transforms via subordination functions, by invoking 
Equation (\ref{eqn:2.23}) from subsection 2.7.

Let $F_{\mathrm{lhs}}$ denote the reciprocal Cauchy transform of the 
probability measure on the left-hand side of (\ref{eqn:3.1}).
We have:
\begin{align*}
F_{\mathrm{lhs}} (z) 
 &= (1-q)z + q F_{\mu^{\boxplus p}} (z) 
      & \mbox{  (by Equation (\ref{eqn:2.170}))}                  \\ 
 &= (1-q)z + q \cdot \frac{p \cdot \omega_{\mathrm{lhs}} (z) -z}{p-1} 
      & \mbox{  (by Equation (\ref{eqn:2.23}))}                   \\ 
 &= \Bigl( (1-q) - \frac{q}{p-1} \Bigr) z + 
        \frac{pq}{p-1} \omega_{\mathrm{lhs}} (z)      &            \\
 &= - \frac{1}{p' -1} z + \frac{p'}{p'-1} 
                   \omega_{\mathrm{lhs}} (z),         &
\end{align*}
where $\omega_{\mathrm{lhs}}$ denotes the subordination function of 
$\mu^{\boxplus p}$ with respect to $\mu$. At the last equality sign 
in the above calculation we used the fact that, due to how $p'$ 
is defined, we have 
$(1-q) - q/(p-1) = -1/(p'-1)$ and $pq/(p-1) = p'/(p'-1)$.

On the other hand, let $F_{\mathrm{rhs}}$ be the reciprocal Cauchy 
transform of the probability measure on the right-hand side of 
(\ref{eqn:3.1}). If we also make the notation
$\nu := \mu^{\uplus q'}$,
then the formula (\ref{eqn:2.23}) gives us that
\[
F_{\mathrm{rhs}} (z) = 
\frac{ p' \cdot \omega_{\mathrm{rhs}} (z) - z}{p'-1}
= - \frac{1}{p'-1} z + \frac{p'}{p'-1} \omega_{\mathrm{rhs}} (z),
\]
where $\omega_{\mathrm{rhs}}$ is the subordination function of 
$\nu^{\boxplus p'}$ with respect to $\nu$. By comparing the latter
expression with the one obtained in the preceding 
paragraph, we see that the desired equality 
$F_{\mathrm{lhs}} = F_{\mathrm{rhs}}$ is tantamount to the equality of 
subordination functions 
$\omega_{\mathrm{lhs}} = \omega_{\mathrm{rhs}}$.

Now, in order to prove that 
$\omega_{\mathrm{lhs}} = \omega_{\mathrm{rhs}}$ we invoke the 
equivalence (\ref{eqn:2.26}) from subsection 2.7,
which tells us that it suffices to check the equality of the left 
inverses of these subordination functions. Let us denote 
these inverses by $H_{\mathrm{lhs}}$ and $H_{\mathrm{rhs}}$, 
respectively. From Equation (\ref{eqn:2.25}) we know that 
\begin{equation}   \label{eqn:3.5}
H_{\mathrm{lhs}}  (w) = pw   + (1-p)  F_{\mu} (w), \ \ 
H_{\mathrm{rhs}} (w) = p' w + (1-p') F_{\nu} (w),
\end{equation}
and the equality $H_{\mathrm{lhs}} = H_{\mathrm{rhs}}$ thus amounts 
to the fact that
\begin{equation}   \label{eqn:3.6}
F_{\nu} (w) = \frac{(p-p')w + (1-p)  F_{\mu} (w)}{1-p'} .
\end{equation}
But on the other hand, since $\nu$ is defined as $\mu^{\uplus q'}$, 
formula (\ref{eqn:2.170}) tells us that 
\begin{equation}  \label{eqn:3.65}
F_{\nu} (w) = (1-q') w + q' F_{\mu} (w);
\end{equation}
and it is straightforward to check that (\ref{eqn:3.6}) reduces to 
(\ref{eqn:3.65}) -- that is, we have
$(p-p')/(1-p') = 1-q'$ and $(1-p)/(1-p') = q'$, due to how $p'$ and
$q'$ were defined in terms of $p$ and $q$.
\hfill$\square$
\end{proof}

\begin{remark} \label{rem:3.20}
{\rm For a measure $\mu \in \cM$ which has compact support, the above 
proposition can be also proved by using combinatorial methods, on the line
developed in \cite{etac} (which has the merit that they extend to a 
multi-variable framework - see Proposition 4.2 in \cite{IJM}).
Here we preferred the treatment via analytic 
methods, which apply directly to a general measure $\mu \in \cM$. }
\end{remark}

$\ $

We now turn to examine why every transformation $\bB_t$ is a 
homomorphism with respect to free multiplicative convolution.
We will prove that, in fact, each of the two kinds of convolution 
powers involved in the definition of $\bB_t$ is ``only a dilation 
away'' from being a homomorphism with respect to $\boxtimes$. In 
order to do that, it will be convenient to start by recording how 
the various transforms considered in the paper behave under dilation.

\begin{notation}  \label{def:3.2}
For $\mu \in {\mathcal M}$ and $r>0$ we denote by $\mu \circ D_r$ 
the probability measure on ${\mathbb R}$ defined by
\[
( \mu \circ D_r ) (A) := \mu (rA), \ \ \forall \, A \subseteq 
{\mathbb R}, \ \mbox{Borel set.}
\]
\end{notation}

\begin{remark}  \label{rem:3.3}
{\rm Let $\mu$ be a probability measure in $\cM$ and let $r$ be a 
positive real number. Then the Cauchy transform of the dilated measure 
$\mu \circ D_r$ is given by the formula
\begin{equation}  \label{eqn:3.7}
G_{\mu\circ D_r}(z) = r G_\mu(rz), \ \ z\in \bC^+.
\end{equation}
This follows directly from the definition of the Cauchy transform:
\[
G_{\mu\circ D_r}(z) = \int_{\bR} \frac{d\mu\circ D_r(x)}{z-x}
= \int_{\bR} \frac{d\mu(x)}{z-\frac{x}{r}}
= r \int_{\bR} \frac{d\mu(x)}{rz-x}
= r G_\mu(rz).
\]
${ }$From Equation (\ref{eqn:3.7}) one easily obtains formulas for 
how various other transforms considered in this paper
change under dilations. For the $R$-transform and the $S$-transform
it turns out that we have
\begin{equation}  \label{eqn:3.71}
R_{\mu \circ D_r} (z) = R_{\mu} (z/r), \ \ 
S_{\mu \circ D_r} (z) = r S_{\mu} (z);
\end{equation}
while for the $\eta$-transform and the $\Sigma$-transform we have 
\begin{equation}  \label{eqn:3.72}
\eta_{\mu \circ D_r} (z) = \eta_{\mu} (z/r), \ \ 
\Sigma_{\mu \circ D_r} (z) = r \Sigma_{\mu} (z).
\end{equation}
Each of the formulas listed in (\ref{eqn:3.71}) and 
(\ref{eqn:3.72}) holds for $z$ running in 
the appropriate domain where the corresponding transform is defined.
The verification of these formulas is immediate, and left as 
exercise to the reader (one just has to start from (\ref{eqn:3.7}) 
and then move through the definitions of the other four transforms 
appearing in the formulas). }
\end{remark}

\begin{proposition}  \label{prop:3.4} 
For $\mu \in {\mathcal M}$ and $\nu \in {\mathcal M}_{+}$ we have
\begin{equation}  \label{eqn:3.81}
( \mu^{\boxplus t} ) \boxtimes ( \nu^{\boxplus t} )  = 
( \mu \boxtimes \nu )^{\boxplus t} \circ D_{1/t},
\ \   \forall \, t \geq 1,                    
\end{equation}
and
\begin{equation}  \label{eqn:3.82}
( \mu^{\uplus t} ) \boxtimes ( \nu^{\uplus t} ) = 
( \mu \boxtimes \nu )^{\uplus t} \circ D_{1/t}, \ \ \forall \, t > 0.
\end{equation}       
\end{proposition}

\begin{proof}
By using the continuity of the operations $\boxplus,\uplus,$ and 
$\boxtimes$ with respect to the weak topology, and by doing suitable
approximations of $\mu,\nu$ in this topology, we see that it suffices
to verify the required formulas in the case when $\mu$ and 
$\nu$ have compact support and have first moment different from 0.
We fix for the whole proof two such measures $\mu$ and $\nu$. For 
these $\mu$ and $\nu$ we will verify (\ref{eqn:3.81}) and 
(\ref{eqn:3.82}) by using the $S$-transform and respectively the 
$\Sigma$-transform.

In order to start on the proof of (\ref{eqn:3.81}), it is useful to 
record how the $S$-transform behaves under $\boxplus$-convolution 
powers: if $\rho \in \cM$ has compact support and non-vanishing first
moment, then we have 
\begin{equation}  \label{eqn:3.9}
S_{\rho^{\boxplus t}} (z) = \frac{1}{t} S_{\rho} (z/t), 
\ \ \forall \, t \geq 1.
\end{equation}
This formula follows immediately from how the $S$-transform is 
expressed in terms of the $R$-transform (that is,
$S (z) = \frac{1}{z} R^{-1} (z)$), combined with the fact 
that $R_{\rho^{\boxplus t}} = t R_{\rho}$.

Now, let us verify that the probability measures appearing on the two
sides of Equation (\ref{eqn:3.81}) have indeed the same $S$-transform.
We compute:
\begin{align*}
S_{( \mu^{\boxplus t} ) \boxtimes ( \nu^{\boxplus t} ) } (z) 
& = S_{\mu^{\boxplus t}} (z) \cdot S_{\nu^{\boxplus t}} (z) \ 
    \mbox{ (by multiplicativity of $S$-transform)}             \\
& = \frac{1}{t^2} S_{\mu} (z/t) \cdot S_{\nu} (z/t) \ 
       \mbox{ (by Equation (\ref{eqn:3.9})), }
\end{align*}
and on the other hand
\begin{align*}
S_{(\mu \boxtimes \nu)^{\boxplus t} \circ D_{1/t}} (z) 
& = \frac{1}{t} S_{(\mu \boxtimes \nu)^{\boxplus t}} (z) 
\ \ \mbox{ (by the $S$-transform formula in 
            Equation (\ref{eqn:3.71})) }                        \\
& = \frac{1}{t^2} S_{\mu \boxtimes \nu} (z/t)
     \ \ \mbox{ (by Equation (\ref{eqn:3.9}))}                  \\
& = \frac{1}{t^2} S_{\mu} (z/t) \cdot S_{\nu} (z/t)
     \ \mbox{ (by multiplicativity of $S$-transform).}
\end{align*}
This completes the verification of (\ref{eqn:3.81}).

The proof of Equation (\ref{eqn:3.82}) goes on the same lines 
as above, with the difference that we now use $\eta$ and $\Sigma$ 
instead of $R$ and $S$. We first note the analogue of Equation 
(\ref{eqn:3.9}):
\begin{equation}  \label{eqn:3.10}
\Sigma_{\rho^{\uplus t}} (z) = \frac{1}{t} \Sigma_{\rho} (z/t), 
\ \ \forall \, t > 0,
\end{equation}
holding for same kind of $\rho$ as in (\ref{eqn:3.9}). Formula 
(\ref{eqn:3.10}) follows immediately from how the $\Sigma$-transform 
is expressed in terms of the $\eta$-function (that is,
$\Sigma_{\rho} (z) = \frac{1}{z} \eta^{-1} (z)$), combined with the 
fact that $\eta_{\rho^{\uplus t}} = t \eta_{\rho}$. By using 
(\ref{eqn:3.10}) and the multiplicativity of the $\Sigma$-transform 
with respect to $\boxtimes$, we obtain (by calculations which are 
virtually identical to those shown in the preceding paragraph) that 
the probability measures on both sides of 
(\ref{eqn:3.82}) have the same $\Sigma$-transform, equal to 
$\frac{1}{t^2} \Sigma_{\mu} (z/t) \cdot \Sigma_{\nu} (z/t)$.
\hfill$\square$
\end{proof}

\begin{remark}  \label{rem:3.45}
{\rm It is immediately seen that Equation (\ref{eqn:3.82})
from the above proposition can also be put in the alternative form 
\begin{equation}  \label{eqn:3.820}
\Bigl( \ ( \mu^{\uplus t} ) \boxtimes ( \nu^{\uplus t} ) 
\ \Bigr)^{\uplus 1/t} = 
( \mu \boxtimes \nu ) \circ D_{1/t}, \ \ \forall \, t > 0 
\end{equation}       
(where $\mu \in \cM$ and $\nu \in \cM_+$, as in 
Proposition \ref{prop:3.4}). 
We note here that in the terminology of \cite{BM}, the 
left-hand side of Equation (\ref{eqn:3.820}) would be said to
define ``the $t$-transform of free multiplicative convolution'';
so in this terminology, (\ref{eqn:3.820}) shows that the 
$t$-transform of $\boxtimes$ is simply obtained by dilating 
$\boxtimes$ by a factor of $1/t$.}
\end{remark}

\begin{remark}  \label{rem:3.5}
(Proof of Theorem 1). 
\rm{ 
At this moment it has become quite easy to verify all
the properties of the transformations $\bB_t$ that were stated in
Theorem 1 of the introduction. Indeed, every $\bB_t$ is injective,
continuous and carries $\cM_{+}$ itself because each of the maps
$$\cM \ni \mu \mapsto \mu^{\boxplus (1+t)} \in \cM \mbox{ and }
\cM \ni \mu \mapsto \mu^{\uplus 1/(1+t)} \in \cM$$
has these properties. The fact that $\bB_t$ is a 
homomorphism with respect to $\boxtimes$ is a straightforward 
consequence of Proposition \ref{prop:3.4}: the dilation factors 
which appear when we take succesively the powers ``$\boxplus (t+1)$'' 
and ``$\uplus 1/(t+1)$'' cancel each other, and we are left with the 
plain multiplicativity stated in Equation (\ref{eqn:1.5}). Finally,
the formula $\bB_t \circ \bB_s = \bB_{t+s}$ 
follows from a direct application of Proposition \ref{prop:3.1}:
\begin{eqnarray*}
\bB_t (\bB_s(\mu)) 
& = & \bB_t \left( (\mu^{\boxplus s+1})^{\uplus\frac{1}{s+1}}
            \right)                                                \\
& = & \left[ \left( (\mu^{\boxplus s+1})^{\uplus\frac{1}{s+1}}
             \right)^{\boxplus t+1} \right]^{\uplus\frac{1}{t+1}}  \\
& = & \left[ \left (\mu^{\boxplus s+1})^{\boxplus
       \frac{t+s+1}{s+1}} \right)^{\uplus\frac{t+1}{t+s+1}}
             \right]^{\uplus\frac{1}{t+1}}                         \\
& = & \left( \mu^{\boxplus{t+s+1}}\right)^{\uplus
               \frac{1}{t+s+1}}                                    \\
& = & \bB_{t+s}(\mu),
\end{eqnarray*}
where at the third equality sign we used Proposition \ref{prop:3.1},
with $p = (t+s+1)/(s+1)$ and $q = (t+1)/(t+s+1)$. }
\hfill$\square$
\end{remark}

\begin{remark} \label{rem:3.55}
{\rm As stated in Theorem 1, every $\bB_t$ maps $\cM_+$ into itself; 
but let us point out here that the converse implication
``$\bB_t ( \mu ) \in \cM_+ \Rightarrow \mu \in \cM_+$'' is not 
true in general. Indeed, it may happen that a measure 
$\mu \in \cM \setminus \cM_+$ has $\mu \boxplus \mu \in \cM_+$.
It is true on the other hand that convolution powers with respect 
to $\uplus$ always preserve $\cM_+$ (that is, 
$\nu^{\uplus t} \in \cM_+$ for every $\nu \in \cM_+$ and every 
$t > 0$). So if $\mu \boxplus \mu \in \cM_+$ then it follows that 
\[
\bB_1 ( \mu ) = ( \mu \boxplus \mu )^{\uplus 1/2} \in \cM_+ ,
\]
even though $\mu$ itself might not belong to $\cM_+$.

In order to give a concrete example of probability measure 
$\mu \in \cM \setminus \cM_+$ such that $\mu \boxplus \mu \in \cM_+$,
one can take for instance a suitable $\boxplus$-power of the 
non-symmetric Bernoulli distribution 
$\mu_o = \frac{1}{4} \delta_{-1} + \frac{3}{4} \delta_1$. We leave 
it as an exercise to the reader to check (by computing explicitly 
the necessary Cauchy transforms and $R$-transforms) that 
$\mu_o^{\boxplus t}$ is in $\cM_+$ for sufficiently large $t>1$.
The set $\{ n \in {\mathbb N} \mid n \geq 1,$
$\mu_o^{\boxplus 2^n} \in \cM_+ \}$ is therefore non-empty. This 
set has a minimal element $m$, and the probability measure 
$\mu := \mu_o^{\boxplus 2^{m-1}}$ has the property that 
$\mu \not\in \cM_+$, but $\mu \boxplus \mu \in \cM_+$. (The explicit
calculations of transforms that were left as exercise can in fact 
be pursued to yield the explicit value for $m$ -- one has $m=4$, 
thus the measure $\mu$ of this example is $\mu_o^{\boxplus 8}$.) }
\end{remark}

\begin{remark} \label{rem:3.6}
{\rm The Equations (\ref{eqn:3.9}) and (\ref{eqn:3.10}) shown during
the proof of Proposition \ref{prop:3.4} can be used in order to obtain
a formula which expresses directly the $S$-ransform of $\bB_t ( \mu )$
in terms of the $S$-transform of $\mu$.

More precisely, let $\mu$ be a probability measure in $\cM$ which has
compact support and has first moment different from 0. Then for every
$t\geq 0$, the measure $\bB_t ( \mu )$ also has these two properties
(this happens because taking convolution powers with respect to 
$\boxplus$ or to $\uplus$ preserves the compactness of the support, 
and rescales the first moment by a factor equal to the exponent used 
in the convolution power). Thus, for every $t \geq 0$, it makes 
sense to consider the $S$-transform of $\bB_t ( \mu )$. We claim that 
for $z$ in a sufficiently small disc centered at 0, this $S$-transform
is given by the formula:
\begin{equation}  \label{eqn:3.11}
S_{\bB_t(\mu)} (z) = S_{\mu} \Bigl( \frac{z}{1-tz} \Bigr) .
\end{equation}  
Indeed, for $|z|$ small enough we can write:
\begin{align*}
S_{\bB_t ( \mu )} (z) 
& = S_{ (      \mu^{ \boxplus (1+t)} )^{\uplus 1/(1+t)} } (z)      \\
& = \Sigma_{ ( \mu^{ \boxplus (1+t)} )^{\uplus 1/(1+t)} }
  ( \frac{z}{1+z} ) \ \mbox{ (by relation between $S$ and $\Sigma$, 
                              Eqn.(\ref{eqn:2.140}))}                \\
& = (1+t) \cdot \Sigma_{ \mu^{ \boxplus (1+t)} } \Bigl( (1+t) \cdot 
  \frac{z}{1+z} \Bigr) \ \mbox{ (by Equation (\ref{eqn:3.10})}       \\
& = (1+t) \cdot S_{ \mu^{ \boxplus (1+t)} } \Bigl( 
  \frac{(1+t) \cdot z/(1+z)}{1- (1+t) \cdot z/(1+z)} \Bigr) \
  \mbox{ (by relation between $S$ and $\Sigma$)}                    \\
& = (1+t) \cdot S_{ \mu^{ \boxplus (1+t)}  } \Bigl( 
    (1+t) \cdot \frac{z}{1-tz} \Bigr)                              \\
& = S_{\mu} \Bigl( \frac{z}{1-tz} \Bigr) 
    \mbox{ (by Equation (\ref{eqn:3.9})). }       
\end{align*}

Note that the multiplicativity of $\bB_t$ with respect to 
$\boxtimes$
could also be inferred from Equation (\ref{eqn:3.11}) (as an 
alternative to the more detailed formulas presented in Proposition 
\ref{prop:3.4}).  }
\end{remark}

\begin{remark}  \label{rem:3.7}
(Proof of Theorem 2).
{\rm At this moment it has become immediate to verify
the statement of Theorem 2 from the introduction. Indeed, since
both $\bB$ and $\bB_1$ are continuous with respect to 
weak topology on $\cM$, it will suffice to verify the equality 
$\bB ( \mu ) = \bB_1 ( \mu )$ in the case when $\mu$ has compact 
support and has first moment different from 0. So let us fix 
$\mu \in \cM$ with these properties. As reviewed in Equation 
(\ref{eqn:2.21}) of subsection 2.6, we have that
\[
S_{\bB ( \mu )} (z) = \Sigma_{\mu} (z), \ \ 
\mbox{ for $|z|$ small enough.} 
\]
But on the other hand, for $|z|$ small enough we also have that
\begin{align*} 
S_{\bB_1 ( \mu )} (z) 
& = S_{\mu} \Bigl( \frac{z}{1-z} \Bigr) \ \mbox{ (by Equation
       (\ref{eqn:3.11}) in Remark \ref{rem:3.6}) }              \\
& = \Sigma_{\mu} (z) \ \mbox{ (by the relation between $S$ and 
                               $\Sigma$, Equation (\ref{eqn:2.140})).}
\end{align*}
Thus the measures $\bB ( \mu )$ and $\bB_1 ( \mu )$ must indeed be 
equal to each other, since they have the same $S$-transform.
\hfill$\square$  }
\end{remark}

\section{Relation to free Brownian motion and to complex Burgers 
         equation}

We first record, in the following lemma, an analytic description 
\setcounter{equation}{0}
of $\bB_t$.

\begin{lemma} \label{b}
Consider $\mu \in\cM$ and $t\geq0$, and consider the measure 
$\bB_t ( \mu ) \in \cM$. We have
\begin{equation}  \label{eqn:4.1}
F_{\bB_t(\mu)} (z) =
\left(1-\frac{1}{t}\right) z + \frac{1}{t} \omega (z),
\quad z\in\mathbb C^+,
\end{equation}
where $\omega$ is the subordination function of 
$\mu^{\boxplus (1+t)}$ with respect to $\mu$.
\end{lemma}

\begin{proof}
Since $\bB_t ( \mu )$ = 
$\bigl( \mu^{\boxplus (t+1)} \bigr)^{\uplus 1/(t+1)}$, Equation
(\ref{eqn:2.170}) from subsection 2.4 gives us that
\begin{equation}  \label{eqn:4.2}
F_{\bB_t (\mu)} (z) =
\left( 1-\frac{1}{1+t} \right) z +
\frac{1}{1+t} F_{\mu^{\boxplus (1+t) }} (z),
\quad z\in\mathbb C^+.
\end{equation}
On the other hand, Equation (\ref{eqn:2.23}) from subsection 2.7 
(with $t$ replaced by $t+1$) tells us that 
\begin{equation}  \label{eqn:4.3}
F_{\mu^{\boxplus (t+1)}} (z) = 
\frac{ (t+1) \omega (z) - z }{t}, \ \ z \in \bC^+.
\end{equation}
Substituting (\ref{eqn:4.3}) into (\ref{eqn:4.2}) leads to the 
required formula (\ref{eqn:4.1}).
\hfill$\square$
\end{proof}

\begin{remark}  \label{rem:4.2}
(Proof of Theorem 3).
{\rm Throughout this proof it will come in handy to use the notation 
$h_{\nu} (z) := F_{\nu} (z) -z$, for $\nu \in \cM$ and $z \in \bC^+$.

Let $\mu$ and $h$ be as in the statement of Theorem 3.
For every $t \geq 1$, let $\omega_t$ denote the subordination function
of $\mu^{\boxplus t}$ with respect to $\mu$, and let $H_t$ be the 
left-inverse for $\omega_t$ defined as in Equation (\ref{eqn:2.25}) from
Section 2. Let us observe that, by
Lemma \ref{b} and the definition of $h_\mu$, we have that
$\omega_{t+1}(z)=z+th_\mu(\omega_{t+1}(z))$,
for all $z\in\mathbb C^+$, $t\ge0.$ Thus, using again Lemma \ref{b},
we find that:
\begin{eqnarray}\label{5.26}
h(t,z)
& = & h_{\id_t(\mu)}(z) \nonumber \\
& = & F_{\id_t(\mu)}(z)-z\nonumber\\
& = & \frac1t(\omega_{t+1}(z)-z)\nonumber\\
& = & h_\mu(\omega_{t+1}(z)),\quad z\in
\mathbb C^+, t\ge0.
\end{eqnarray}
We observe that the function $\omega_{t+1}(z)$ is indeed differentiable 
in both variables; this follows immediately from the equation 
$H_{t+1}(\omega_{t+1}(z))=z$ and the definition of $H_{t+1}(z)$. Let us 
denote $\frac{\partial}{\partial x}$ by $\partial_x$. Differentiating 
with respect to $z$ gives

$$\partial_z\omega_{t+1}(z)=
\frac{1}{(\partial_zH_{t+1})(\omega_{t+1}(z))},$$
and differentiating with respect to $t$ gives

$$\partial_t\omega_{t+1}(z)=
-\frac{(\partial_tH_{t+1})(\omega_{t+1}(z))}{(\partial_zH_{t+1})(\omega_{t+1}(z))}
={h_\mu(\omega_{t+1}(z))}{\partial_z\omega_{t+1}(z)},$$
where we have used the definitions of $H$ and $h$ in the last equality.

Then
$$\partial_th(t,z)=\partial_th_{\id_t(\mu)}(z)=
h_\mu'(\omega_{t+1}(z))\partial_t\omega_{t+1}(z)=
h_\mu'(\omega_{t+1}(z)){h_\mu(\omega_{t+1}(z))}{\partial_z\omega_{t+1}(z)}
,$$
and
$$\partial_zh(t,z)=\partial_zh_{\id_t(\mu)}(z)=
h_\mu'(\omega_{t+1}(z))\partial_z\omega_{t+1}(z)
.$$
The two relations above prove equation (\ref{eqn:1.9}).
\hfill$\square$  }
\end{remark}

\begin{remark}  \label{rem:4.3}
{\rm We now move to the framework of Theorem 4. Let $\nu$ be in $\cM$,
and let us consider the Cauchy transform $G_{\nu}$. It is easily 
verified that the map $\bC^+ \ni z \mapsto z - G_{\nu} (z) \in \bC$ 
belongs to the set $\mathcal F$ of analytic self-maps of $\bC^+$ 
considered in the above Equation (\ref{eqn:2.3}); thus there exists 
a unique $\mu \in \cM$ such that $z - G_{\nu} (z) = F_{\mu} (z)$, 
$z \in \bC^+$. Clearly, this $\mu$ is related to the given $\nu$ in 
exactly the way described by Equation (\ref{eqn:1.12}) of Theorem 4. 
It can be shown that $\mu$ has variance equal to 1 and is
centered, that is, it satisfies
\begin{equation}   \label{eqn:4.10}
\int_{- \infty}^{\infty} t^2 \ d \mu (t) = 1, \ \ 
\int_{- \infty}^{\infty} t \ d \mu (t) = 0. 
\end{equation}
Moreover, it can be shown that the correspondence $\nu \mapsto \mu$ 
(with $\nu$ and $\mu$ as above) is a bijection between $\nu \in \cM$ 
and $\mu$ running in the set of probability measures in $\cM$ which 
satisfy the conditions in (\ref{eqn:4.10}). A detailed presentation 
of these facts appears in Section 2 of the paper \cite{M92} of 
Maassen (see Proposition 2.2 of that paper).  }
\end{remark}

\begin{remark}  \label{rem:4.4}
(Proof of Theorem 4).
{\rm Let us fix two probability measures $\mu , \nu \in \cM$ which 
are connected to each other as in Remark \ref{rem:4.3}, via the 
relation
\[
-G_{\nu} (z) = F_{\mu} (z) -z, \ \ z \in \bC^+ .
\]
Let us also fix a real number $t > 0$. Our goal in this proof is to
show that we have the formula
\begin{equation}   \label{eqn:4.11}
-G_{\nu \boxplus \gamma_t} (z) = F_{\bB_t ( \mu )} (z) -z, 
\ \ z \in \bC^+ ,
\end{equation}
where $\gamma_t$ is the centered semicircular distribution of 
variance $t$. 

Let $\theta : \bC^+ \to \bC^+$ be the subordination function of 
$\nu \boxplus \gamma_t$ with respect to $\nu$. By the definition 
of $\theta$, we thus have 
\begin{equation}   \label{eqn:4.12}
G_{\nu \boxplus \gamma_t}  = G_{\nu} \circ \theta
\end{equation}
(composition of self-maps of $\bC^+$). But, as noted at the end of
the above subsection 2.7, the maps $\theta$ and 
$G_{\nu \boxplus \gamma_t}$ are also related via the formula
\begin{equation}   \label{eqn:4.13}
\theta (z) = z - t G_{\nu \boxplus \gamma_t} (z), \ \ z \in \bC^+ . 
\end{equation}

Now let $\omega : \bC^+ \to \bC^+$ be the subordination function 
of $\mu^{\boxplus (t+1)}$ with respect to $\mu$. We claim that 
$\omega$ is in fact equal to the subordination function $\theta$ from 
the preceding paragraph. In order to prove this claim, we write the 
Cauchy transform of $\nu \boxplus \gamma_t$ in two different ways.
On one hand, for every $z \in \bC^+$ we have
$G_{\nu \boxplus \gamma_t} (z) = G_{\nu} ( \theta (z) )$
(by Equation (\ref{eqn:4.12})), hence
\begin{equation}   \label{eqn:4.14}
G_{\nu \boxplus \gamma_t} (z) = - F( \, \theta (z) \, ) + \theta (z)
\end{equation}
(where we used the formula which connects $G_{\nu}$ to $F_{\mu}$, 
applied to the complex number $\theta (z) \in \bC^+$). On the other 
hand, from Equation (\ref{eqn:4.13}) we have that
\begin{equation}   \label{eqn:4.15}
G_{\nu \boxplus \gamma_t} (z) = \frac{1}{t} \Bigl( 
z - \theta (z) \Bigr) .
\end{equation}
By eliminating $G_{\nu \boxplus \gamma_t} (z)$ between Equations
(\ref{eqn:4.14}) and (\ref{eqn:4.15}), we find that
\[
- F( \, \theta (z) \, ) + \theta (z) =
\frac{1}{t} \Bigl( z - \theta (z) \Bigr) ,
\]
which implies that
\begin{equation}   \label{eqn:4.16}
\theta (z) = \frac{1}{t+1} z + \Bigl( 1 - \frac{1}{t+1} \Bigr)
F_{\mu} ( \, \theta (z) \, ), \ \ z \in \bC^+ .
\end{equation}
We obtained that $\theta$ satisfies the functional 
equation which determines uniquely the subordination function 
$\omega$ (cf. discussion in subsection 2.7 above); the equality 
$\theta = \omega$ follows.

Finally, let us observe that for every $z \in \bC^+$ we have
\begin{align*}
F_{\bB_t ( \mu )} (z) - z 
& = \Bigl( ( 1 - \frac{1}{t} ) z + \frac{1}{t} \omega (z) 
    \Bigr) -z \ \mbox{ (by Lemma \ref{b})}                    \\
& = \frac{1}{t} \Bigl( \omega (z) - z \Bigr)                  \\
& = \frac{1}{t} \Bigl( \theta (z) - z \Bigr) \ \mbox{ (since
           $\omega = \theta$) }                               \\
& = 
- G_{\nu \boxplus \gamma_t} (z) \ \mbox{ (by Equation 
                                         (\ref{eqn:4.15})), }  
\end{align*}
and the desired Equation (\ref{eqn:4.11}) is obtained.
\hfill$\square$    }
\end{remark}

\begin{example}  \label{ex:4.5}
{\rm  Let $\mu$ be the symmetric Bernoulli distribution,
$\mu = \frac{1}{2} ( \delta_{-1} + \delta_1 )$. The measure 
$\nu$ corresponding to this $\mu$ via the bijection 
$\nu \leftrightarrow \mu$ from Remark \ref{rem:4.3} is the 
Dirac measure $\delta_0$; indeed, with $\nu = \delta_0$,
we clearly have $- G_{\nu} (z) = F_{\mu} (z) -z = -1/z$,
$z \in \bC_+$. But then Theorem 4 implies that for every 
$t \geq 0$ we have 
\begin{align*}
F_{ \bB_t ( \mu ) } (z) - z 
& = - G_{ \delta_0 \boxplus \gamma_t } (z)              \\
& = - G_{ \gamma_t } (z)                                \\
& = -\frac{z-\sqrt{z^2-4t}}{2t}, \ z \in \bC_+.
\end{align*}
We thus get an explicit formula for $F_{ \bB_t ( \mu ) }$, 
leading to a formula for the Cauchy transform of $\bB_t ( \mu )$: 
\begin{equation}  \label{eqn:4.17}
G_{ \bB_t ( \mu ) } (z) = 
\frac{ (2t-1) z - \sqrt{z^2 - 4t} }{ 2(1- (1-t)z^2) },  
         \ \ z \in \bC^+.
\end{equation}
By a straightforward application of the Stieltjes inversion formula, 
one can then determine exactly what $\bB_t ( \mu )$ is. For 
$t \geq 1/2$ the result of the calculation is that $\bB_t ( \mu )$
is absolutely continuous, with density
\begin{equation}  \label{eqn:4.18}
x \mapsto \frac{ \sqrt{4t-x^2} }{ 2 \pi \cdot (1-(1-t)x^2 ) }, 
                            \ \ |x| \leq 2 \sqrt{t}.
\end{equation}
For $0 < t < 1/2$, one finds that $\bB_t ( \mu )$ has an 
absolutely continuous part with density described exactly as 
in the above Equation (\ref{eqn:4.18}); and in addition to that,
$\bB_t ( \mu )$ has two atoms at $\pm 1/\sqrt{1-t}$, each of them of 
mass equal to $(1-2t)/(2-2t)$. (The details of the calculation are 
left as an exercise to the reader.)

Note that for $t = 1/2$ the density from (\ref{eqn:4.18}) 
simplifies to  
\[
x \mapsto \frac{1}{ \pi \sqrt{2-x^2} }, \ \  |x| \leq \sqrt{2}. 
\]
Thus $\bB_{1/2} ( \mu )$ is the arcsine law on the interval 
$[ - \sqrt{2} , \sqrt{2} ]$. 

On the other hand, for $t=1$ the density from 
(\ref{eqn:4.18}) becomes
\[
x \mapsto \frac{ \sqrt{4-x^2} }{ 2 \pi }, \ \  |x| \leq 2, 
\]
which shows that $\bB_{1} ( \mu )$ is the standard semicircle law
$\gamma_1$. (The fact that $\bB_1 ( \mu ) = \gamma_1$ could also 
be obtained directly from Theorem 2. Indeed, $\bB_1 ( \mu )$ is 
equal to $\bB ( \mu )$, and is hence determined by the fact that 
it has $R$-transform equal to $\eta_{\mu} (z)$. But it is immediate
that $\mu$ of this example has $\eta_{\mu} (z) = z^2$, which is 
exactly the $R$-transform of $\gamma_1$.)   }
\end{example}

\begin{remark}   \label{rem:4.6}{\rm
The above Example \ref{ex:4.5} can be placed within the framework 
of a family of probability measures on $\bR$ which are called 
{\em free Meixner states}. (We thank Michael Anshelevich for 
bringing this observation to our attention.) We briefly outline 
here how the connection to free Meixner states appears;
for more
details, see Section 4.2 in the recent paper \cite{A}.

Let $b$ and $c$ be two real parameters such that $c \geq -1$. The 
free Meixner state of mean 0 and variance 1 which is indexed by the 
parameters $b,c$ is a probability measure $\mu_{b,c} \in \cM$ 
that can be defined as follows: if $( P_n )_{n=1}^{\infty}$ is the 
sequence of monic orthogonal polynomials for $\mu_{b,c}$, then 
we have $P_0 (t)=1$, $P_1(t)=t$, and the recurrence
\begin{equation}  \label{eqn:4.140}
\left\{  \begin{array}{cccc}
tP_1(t) & = & P_2(t) + bP_1(t) + P_0(t),    &                    \\
tP_n(t) & = & P_{n+1}(t) + bP_n(t) + (c+1)P_{n-1}(t),
                                             & \forall \, n \geq 2. 
\end{array}  \right.
\end{equation}
For the explanation of why the probability measures $\mu_{b,c}$ are 
called ``free Meixner'', we refer to the paper \cite{An03}. 

Now, it is not hard to obtain an explicit formula for the Cauchy 
transform of $\mu_{b,c}$ -- see the displayed equations in Theorem 4 
of \cite{An03} (the parameters ``$a,t$'' from those equations have 
to be suitably set, in order to match the $b,c$ from the above 
Equation (\ref{eqn:4.140})). This leads to an explicit formula for 
the analytic function $F_{\mu_{b,c}} (z) -z$; it turns out that 
what one gets is exactly the relation
\begin{equation}  \label{eqn:4.150}
F_{\mu_{b,c}} (z) -z = - G_{\gamma_{b,c+1}} (z), \ \ z \in \bC^+,
\end{equation}
where $\gamma_{b,t}$ stands for the semicircular distribution of 
mean $b$ and variance $t$. Or in other words, the correspondence 
$\nu \mapsto \mu$ from Remark \ref{rem:4.3} has the property that
\begin{equation}  \label{eqn:4.160}
\gamma_{b,t} \mapsto \mu_{b,t-1}, \ \ \forall \, b \in \bR 
\mbox{ and } t \geq 0.
\end{equation}
(Thus the set of all free Meixner states with mean 0 and variance 1 
is precisely what one obtains when applying the correspondence from 
Remark \ref{rem:4.3} to the set of all -- not necessarily centered
-- semicircular distributions!)

It is of course clear that we have
\[
\gamma_{b,t} = \delta_b \boxplus \gamma_t , \ \ 
\forall \, b \in \bR , \ \ \forall \, t \geq 0,
\]
where $\gamma_t$ is the centered semicircular distribution of 
variance $t$. Thus Theorem 4 of this paper can be applied in 
conjunction to the above Equation (\ref{eqn:4.160}), and it gives 
us that
\begin{equation}  \label{eqn:4.170}
\bB_t  ( \mu_{b,c} ) = \mu_{b,c+t}, \ \ \forall \, b,c,t \in \bR
\mbox{such that } c \geq -1 \mbox{ and } t \geq 0.
\end{equation}

If we fix a value $b \in \bR$, we hence see that the family 
of probability measures $\{ \mu_{b,c} \mid c \geq -1 \}$ is 
exactly what one obtains by starting with the measure 
$\mu_{b,-1}$ and by letting it evolve under the semigroup of 
transformations $\bB_t$. It is easily computed that 
$\mu_{b,-1}$ is actually a Bernoulli measure,
\begin{equation}  \label{eqn:4.180}
\mu_{b,-1}  = \frac{q}{q-p} \delta_p + 
\frac{p}{p-q} \delta_q , 
\end{equation}
with $p,q$ determined from the equations $p+q =b$, $pq=-1$.
The situation from Example \ref{ex:4.5} is obtained by setting
$b=0$ in this discussion, when $\mu_{b,-1}$ becomes the 
symmetric Bernoulli distribution 
$\frac{1}{2} ( \delta_{-1} + \delta_1 )$.}
\end{remark}

\section{Miscellaneous facts about $\bB_t ( \mu )$}

\setcounter{equation}{0}

We start this section with a proposition discussing atoms and
regularity for a measure $\bB_t ( \mu )$, $t>0$. We show that 
densities of such probability measures tend to be quite smooth, 
while singular continuous parts cannot appear.

\begin{proposition}  \label{prop:5.1}
Let $\mu\in\mathcal M$ and any $t>0$. Then:
\begin{enumerate}
\item The singular continuous part of $\id_t(\mu)$
with respect to the Lebesgue measure on the real line is zero.
\item The probability $\id_t(\mu)$ has at most $[1/t]$
atoms for $t\leq 1$ and at most one atom for $t>1$. The point 
$x\in\mathbb R$ is an atom  of $\id_t(\mu)$ if and only if
\[
\lim_{y\to0}F_\mu((1-t)x+iy)
= -tx\quad{\rm and}\quad \lim_{y\to0}F_\mu'((1-t)x+iy)
= \frac{1+t(1-\id_t(\mu)(\{x\}))}{\id_t(\mu)(\{x\})+t(1-
                                            \id_t(\mu)(\{x\}))}. 
\]
(We have denoted by $[a]$ the largest integer less than or 
equal to $a$.)
\item The absolutely continuous part of $\id_t(\mu)$ with respect 
to the Lebesgue measure is zero if and only if $\mu=\delta_c$ for 
some $c\in\mathbb R$.  Moreover, its density
is analytic wherever positive and finite.
\end{enumerate}
\end{proposition}

\begin{proof}
We shall prove first item 2. 
As known from equation (5.7) in \cite{BBercoviciMult},
if $\lim_{y\to0}F_{\id_t(\mu)}(x+iy)=0$, then we have
\begin{equation}\label{atom}
\id_t(\mu)(\{x\})^{-1}=\lim_{y\to0}\frac{F_{\id_t(\mu)}
(x+iy)}{iy}=\lim_{y\to0}F_{\id_t(\mu)}'
(x+iy),\end{equation}
where we use the convention $1/0=\infty$; conversely, for $x$ 
to be an atom of $\id_t(\mu)$ it is required that
$\lim_{y\to0}F_{\id_t(\mu)}(x+iy)=0$, and then (\ref{atom}) holds.
${ }$From Lemma \ref{b} it follows that this is equivalent to
$$\lim_{y\to0}\omega_{t+1}(x+iy)=(1-t)x$$
$$ \ {\rm and}\  $$  $$\lim_{y\to0}\frac{
\omega_{t+1}(x+iy)-(1-t)x}{iy}=
\lim_{y\to0}\omega_{t+1}'(x+iy)=(1-t)+\frac{t}{\id_t(\mu)(\{x\})}.$$
These statements can be seen to be equivalent to
$\lim_{y\to0}H_{t+1}((1-t)x+iy)=x$ and $\lim_{y\to0}H_{t+1}'((1-t)x+iy)=
\left(1-t+\frac{t}{\id_t(\mu)(\{x\})}\right)^{-1}$, since that 
$H_{t+1}(z)=(t+1)z-tF_\mu(z)$
(See also Proposition 4.7 of \cite{BBercoviciMult}).
Using the definition of $H_{t+1}$ and \cite[(5.7)]{BBercoviciMult},
we obtain
$$\lim_{y\to0}F_\mu((1-t)x+iy)=-tx$$
$${\rm and}$$
$$\lim_{y\to0}F_\mu'((1-t)x+iy)
= \frac{t+1}{t}-\frac{\id_t(\mu)(\{x\})}{t[t+\id_t(\mu)(\{x\})(1-t)]}
= 
\frac{1+t(1-\id_t(\mu)(\{x\}))}{\id_t(\mu)(\{x\})+t(1-\id_t(\mu)(\{x\}))}. 
$$

Assume first that $t<1$. Then the function $F_\mu((1-t)z)+tz$, 
$z\in\mathbb C^+.$
maps the upper half-plane into itself,  since $\Im F_\mu((1-t)z)\geq
(1-t)\Im z,$ so if $1-t>0$ we must have $\Im(F_\mu((1-t)z)+tz)>0$ for all
$z\in\mathbb C^+.$ Moreover, the function $G_\lambda(z)=
1/(F_\mu((1-t)z)-tz)$ is the Cauchy transform of a probability measure 
$\lambda$, and $F_\lambda(z)=1/G_\lambda(z)$
satisfies $$\lim_{y\to0}F_\lambda(x+iy)=0,$$
\begin{eqnarray*}
\lim_{y\to0}F_\lambda'(x+iy) & = & t+(1-t)\lim_{y\to0}F_\mu'((1-t)x+iy)\\
& = & 
t+(1-t)\frac{1+t(1-\id_t(\mu)(\{x\}))}{\id_t(\mu)(\{x\})+t(1-\id_t(\mu)(\{x\}))}\\
& = & \frac{1}{\id_t(\mu)(\{x\})+t(1-\id_t(\mu)(\{x\}))}.
\end{eqnarray*}
This holds for any $x$ so that $\lim_{y\to0}F_{\id_t(\mu)}(x+iy)=0$. Thus, 
any such
point $x$ must be an atom of $\lambda$ and
  $\lambda(\{x\})=\id_t(\mu)(\{x\})+t(1-\id_t(\mu)(\{x\})).$
Thus, $\lim_{y\to0}F_{\id_t(\mu)}(x+iy)=0$ and $0<t<1$ implies
$\lambda(\{x\})\geq t.$ Since the total mass of all atoms of $\lambda$ 
cannot exceed
one, we conclude that there exist at most $[1/t]$ points $x$ where
$\lim_{y\to0}F_{\id_t(\mu)}(x+iy)=0$, and in particular at most $[1/t]$ 
atoms of
$\id_t(\mu)$.

If $t\geq1$, then we have proved in Theorem 2 that $\id_t(\mu)$  is 
infinitely divisible,
so, as observed in \cite{BVReg}, it can have at most one atom. This proves 
item 2.

We shall now show that $\id_t(\mu)$ cannot have a nonzero singular 
continuous part.
Observe that since $\id_t(\mu)$ is $\boxplus$-infinitely divisible for all 
$t\geq1$,
according to Theorem 2, the statement is a trivial consequence of 
Proposition 5.1 in
\cite{BBercoviciMult} whenever $t\geq1.$ Thus we will focus again on
the
case $0<t<1.$
Assume that $\id_t(\mu)$ has a nontrivial singular continuous part. As 
noted in
\cite{BBercoviciMult}, for uncountably many points $x$ in the support of
the singular continuous part of $\id_t(\mu)$, we must have
$$\lim_{y\to0}F_{\id_t(\mu)}(x+iy)=0\quad{\rm 
and}\quad\lim_{y\to0}\frac{F_{\id_t(\mu)}
(x+iy)}{iy}=\lim_{y\to0}F_{\id_t(\mu)}'
(x+iy)=\infty.$$
As shown in the proof of item 2, this implies that there exists a 
probability $\lambda$
that has uncountably many atoms of mass at least $t$, an obvious 
contradiction.
This proves item 1.

We prove next item 3. As observed in \cite{BBercoviciMult},
for any $t>0$, the function $F_{\mu^{\boxplus t+1}}$ extends analytically 
through
any point $x$ where $F_{\mu^{\boxplus t+1}}$ has finite nontangential 
limit with
strictly positive imaginary part,  so, by Lemma \ref{b},
the function $F_{\id_t(\mu)}$ extends analytically through any point $x$
where $F_{\mu^{\boxplus t+1}}>0$ has finite nontangential limit with
strictly positive imaginary part.
Since the density of ${\id_t(\mu)}$ is just $-\pi^{-1}\Im (1/F_{\id_t(\mu)})$, its 
analyticity follows
immediately.

It has been observed in Proposition 4.7 of \cite{BBercoviciMult} that 
$\mu^{\boxplus t+1}$ has a nonzero absolutely continuous part for any 
$t>0$ if and only if $\mu$
is not concentrated in one point. Thus, by Lemma \ref{b}, and the remarks 
above, we
conclude that for any $t>0$, $\id_t(\mu)$ has a nonzero absolutely 
continuous part if and only if $\mu$ is not concentrated in one point. 
This proves item 3.
\hfill$\square$
\end{proof}


\begin{corollary}\label{cor:5.2}
Assume that $\mu\in\mathcal M$ is absolutely continuous with respect to the Lebesgue measure on
the real line, with density $f_0$. Denote by $f_t$ the density with respect to the 
Lebesgue measure of the absolutely continuous part of $\id_t(\mu)$. Then $f_t$ converge to $f_0$
almost everywhere when $t$ tends to zero.
\end{corollary}

\begin{proof}
It is known \cite[Theorem 1.3]{Duren} that any bounded function which 
is analytic in the unit disk has nontangential limit at almost all 
points of the unit circle. Since
 $G_\mu$ maps the upper into the lower half-plane, it can be conjugated
to a self-map of the unit disk by using the map $z\mapsto
\frac{z-i}{z+i}$, and thus has nontangential limit at almost all points $x\in\mathbb R$ (denote the
limit function by $G_\mu(x)$). Also, by the upper half-plane version of
\cite[Theorem 1.2]{Duren},
 $\Im G_\mu(x)=-\pi f_0(x)$ for Lebesgue-almost all $x\in\mathbb R$. Since the same holds when
$\mu$ is replaced by $\id_t(\mu)$ and $f_0$ by $f_t$, it is enough to prove that
$\lim_{t\to0} G_{\id_t(\mu)}(x)=G_\mu(x)$ for all $x\in\mathbb R$, except possibly for a set of
zero Lebesgue measure. The statement of the corollary follows from this limit by taking
the imaginary parts.

Now, from the relation between $G_\mu$ and $h_\mu$ it is immediate that it will be 
enough to show 
that for any $x\in\mathbb R$ with the property that the nontangential limit of $h_\mu$ 
at $x$ belongs to the upper half-plane, we have $\lim_{t\to0}h_{\id_t}(x)=h_\mu(x).$
We observe that for any such $x\in\mathbb R$, the function 
$h_x(z):=h_\mu(x+z),$ $z\in\mathbb C^+$, has nontangential limit at zero belonging to
the upper half-plane, and the function $\gamma(t):=th_{\id_t(\mu)}(x)$ satisfies the
condition $\gamma(t)=th_x(\gamma(t))$.
Indeed, the first statement is just a reformulation of the fact that 
the nontangential limit of $h_\mu$ 
at $x$ belongs to the upper half-plane. For the second observation, note that
by Lemma \ref{b} we have $\omega(z)=z+th_{\id_t(\mu)}(z)$, and by Theorem
4.6 in \cite{BBercoviciMult}, this relation extends by continuity to the real line. Thus, 
relation (\ref{5.26}) together with the above implies that $\gamma(t)=th_x(\gamma(t))$.
These conditions have been proved in Lemma 2.13 in \cite{BG} to imply that
$\lim_{t\to0}h_x(\gamma(t))$ exists and equals the nontangential limit at zero
of $h_x$. But this is equivalent to saying that $\lim_{t\to0}h_{\id_t(\mu)}(x)=h_\mu(x).$\hfill$\square$
\end{proof}

\begin{remark}
{\rm
The proposition above gives a full description of the atoms of 
$\bB_t(\mu)$ in terms of the reciprocal Cauchy transform of $\mu$, 
and a quite strict bound on the possible number of such atoms. 
However, it does not guarantee that the lack of atoms for
$\mu$ must translate into a lack of atoms for $\id_t(\mu)$. 
We will provide an example showing that this situation in fact can 
occur. Consider the probability measure $\mu$ whose reciprocal 
Cauchy transform is given by the formula
\[
F_\mu(z)= z-2+i+\frac{z-i}{z+i},\quad z\in\mathbb C^+\cup\mathbb R.
\]
It is easy to see that $F_\mu$ is in fact a rational transformation 
of the whole complex plane, with one pole at $-i$. We observe that 
$F_\mu(1)=-1$ and $F_\mu(x)\in\mathbb C^+$
for all $x\in\mathbb R\setminus\{1\}$ so $\mu$ has no atoms. 

Let $x=2$, $t=1/2.$ Then
$$F_\mu\left(\left(1-\frac12\right)\cdot2\right)
=F_\mu(1)=-1=-\frac12\cdot2,$$
and since $F_\mu'(z)=1+\frac{2i}{(z+i)^2},$ we have
$$F_\mu'\left(\left(1-\frac12\right)\cdot2\right)=F_\mu'(1)=2.$$
According to the proposition above, $\id_\frac12(\mu)(\{2\})=1/3.$  }
\end{remark}

The remaining part of this section is about the
$\boxplus$-divisibility indicator $\phi ( \mu )$. 

\begin{proposition}  \label{prop:5.4}
Let $\mu$ be in $\cM$. We have that $\mu$ is infinitely divisible 
with respect to $\boxplus$ if and only if $\phi ( \mu ) \geq 1$.
\end{proposition}
\begin{proof}
We have already seen that if $\mu$ is infinitely divisible, then $\mu\in\id(\mathcal M)$,
so that $\phi(\mu)\ge1$. Assume now that $\phi(\mu)\ge1$. We shall prove that $\mu$
is infinitely divisible. This is clearly true if $\phi(\mu)>1$, so we may assume without loss of generality that $\phi(\mu)=1$.
Observe that $\lim_{t\downarrow0}\id_t(\mu)=\mu$ in the weak topology, and (by the
definition of $\phi$) $\phi(\id_t(\mu))>1$ for all $t>0$, so that by Theorem 2, 
$\id_t(\mu)$ is $\boxplus$-infinitely divisible. Thus, the statement
of our proposition will be proved once we show that the set of
$\boxplus$-infinitely divisible probability measures is closed in the
topology of weak convergence.

As shown in \cite{BVIUMJ}, Proposition 5.7, weak convergence of probability measures
translates into convergence on compact subsets of the corresponding $R$-transforms:
there exists a cone $\Gamma$ with vertex at zero so that $\mathcal R_\mu$ and
$\mathcal R_{\id_t(\mu)}$ are all defined on $\Gamma$ for $t>0$ small enough and
$\lim_{t\downarrow0}\mathcal R_{\id_t(\mu)}(z)=\mathcal R_{\mu}(z)$ uniformly on compact subsets of
$\Gamma$. But by Theorem 5.10 in the same \cite{BVIUMJ}, $\mu$ is infinitely
divisible if and only if $\mathcal R_{\mu}$ has an analytic extension to the upper half-plane.
As the $R$-transform $\mathcal R$ maps points from $\mathbb C^+$ into
$\mathbb C^+\cup\mathbb R$, the family 
$\{\mathcal R_{\id_t(\mu)}\colon t\ge0\}$ is normal, so, by the classical theorem of 
Montel,  there exists a subsequence 
$t_n\to0$ so that $\mathcal R_{\id_{t_n}(\mu)}$ is convergent on $\mathbb C^+.$
The existence of the limit on $\Gamma$ together with the uniqueness of analytic 
continuation guarantees that $\mathcal R_\mu$ has an analytic extension to the
upper half-plane, and thus $\mu$ is $\boxplus$-infinitely divisible. 
\hfill$\square$
\end{proof}

\begin{remark}   \label{rem:5.5}
(Some basic properties of $\phi ( \mu )$.) 

{\rm 
$1^o$ A useful formula for deriving explicit values of $\phi(\mu)$ is
\begin{equation}\label{eqn:5.1001}
\phi ( \bB_t ( \mu ) ) = t+ \phi ( \mu ),\quad \forall\mu\in\mathcal M,\forall t\in[0,+\infty).
\end{equation}
Indeed, if $\phi(\mu)=0$, then this is obvious from the injectivity of $\id_t$ and the semigroup property.
If $\phi(\mu)>0$, then 
consider a sequence $q_n$ that increases to $\phi(\mu)$ and define $\nu_n$ by
the equation $\id_{q_n}(\nu_n)=\mu$. Clearly
$$\phi(\id_t(\mu))=\phi(\id_{q_n+t}(\nu_n))\ge q_n+t\to t+\phi(\mu)$$
as $n\to\infty$. Thus, $\phi(\id_t(\mu))\ge t+\phi(\mu)$. We claim that the inequality 
cannot be strict. Indeed, if this were not the case, we would find an $\varepsilon$ so that
$\phi(\id_t(\mu))-( t+\phi(\mu))>\varepsilon>0$ and  a measure $\nu_\varepsilon$ so that
$\id_t(\mu)=\id_{t+\phi(\mu)+\varepsilon}(\nu_\varepsilon)$, and thus, by injectivity of $
\id_t$ together with the semigroup property, $\mu=\id_{\phi(\mu)+\varepsilon}(\nu_\varepsilon)$, so that $\phi(\mu)\ge
\phi(\mu)+\varepsilon,$ an obvious contradiction.

$2^o$ Let $\mu \in \cM$ be such that $\phi ( \mu ) =:p > 0$. From 
the definition of $\phi ( \mu )$ and the semigroup property of the 
transformations $\bB_t$ it is immediate that for every 
$0 < q < p$ one can find a probability measure $\nu_q \in \cM$ 
such that $\bB_q ( \nu_q ) = \mu$. However,
as an immediate consequence of $1^o$ above and of Proposition \ref{prop:5.4}, we observe that
a stronger result holds:
\begin{equation}  \label{eqn:5.101}
\exists \, \nu \in \cM \mbox{ such that } \bB_p ( \nu ) = \mu .
\end{equation}
The proof is as follows: assume first that $\phi(\mu)=1$. Then $\mu$
is $\boxplus$-infinitely divisible by Proposition \ref{prop:5.4}, so
$\mu\in\id_1(\mathcal M),$ by Theorem 2, and (\ref{eqn:5.101}) follows.
 If $\phi(\mu)<1$, consider
$\nu=\id_{1-\phi(\mu)}(\mu)$. By $1^o$ above, $\phi(\nu)=1$, so there exists
$\nu_0\in\mathcal M$ so that $\nu=\id_1(\nu_0)$. Thus, $\id_1(\mu)=\id_{\phi(\mu)}(\nu)
=\id_{1+\phi(\mu)}(\nu_0)$, and, from injectivity of $\id_1$ and the semigroup property,
$\mu=\id_{\phi(\mu)}(\nu_0)$, proving our statement for $\phi(\mu)<1$. 
The similar (and easier) case $\phi(\mu)>1$ is left to the reader.

$3^o$ Let $\mu \in \cM$ be such that $0 < \phi ( \mu ) < 1$. Then 
we can take  convolution powers $\mu^{\boxplus t}$ for any
$1-\phi(\mu)\le t\le1$.
Indeed, let us fix such a $t\in[1-\phi(\mu),1]$. Then $0\leq1-t\leq\phi(\mu)$ and hence 
there exists $\nu\in\mathcal M$ so that $\id_{1-t}(\nu)=\mu$.
Now,

\begin{eqnarray*}
\mu & = & \id_{1-t}(\nu)\\
& = & \left(\nu^{\boxplus(2-t)}\right)^{\uplus1/(2-t)}\\
& = & \left(\nu^{\uplus t}\right)^{\boxplus1/t} {\quad\ \ \ \textrm{(by 
Proposition \ref{prop:3.1}}).}
\end{eqnarray*}
Since $t\le1$, the last term in the above equalities
is defined, and thus $\mu^{\boxplus t}$ is well defined by the above and 
equals $\nu^{\uplus t}$.

$4^o$ We have
\begin{equation}   \label{eqn:5.102}
\phi ( \mu \boxtimes \nu ) \geq 
\min \{ \phi ( \mu ), \phi ( \nu ) \}, 
\ \ \forall \, \mu , \nu \in \cM.
\end{equation}
This follows immediately from the homomorphism property of $\id_t$ proved 
in Theorem 1, and from $2^o$ above.
An immediate consequence of this fact is that if both $\mu$ and 
$\nu$ are infinitely divisible with respect to $\boxplus$, then 
so is $\mu \boxtimes \nu$. (This statement is not obvious from the 
definitions, but can also be derived directly from 
Equation (\ref{eqn:3.81}) in Proposition \ref{prop:3.4}, 
without resorting to $\boxplus$-divisibility indicators.)  }
\end{remark}

\begin{remark}   \label{rem:5.6}
{\rm It is now easy to derive the concrete values $\phi ( \mu )$ that were
listed in Table 1 of the introduction section. Indeed, from Proposition 
\ref{prop:5.1} it is immediate that $B_t ( \mu )$ can never have finite 
support when $t > 0$, and this implies in particular that the symmetric 
Bernoulli distribution $\mu_o = \frac{1}{2} ( \delta_1 + \delta_{-1})$
has $\phi ( \mu_o ) = 0$. But for this $\mu_o$ we know that 
(as observed in the above Example \ref{ex:4.5}) the measure 
$\bB_{1/2} ( \mu_o )$ is the arcsine law, while $\bB_1 ( \mu_o )$ is 
the standard semicircle law. The values of $\phi ( \mu )$ listed 
on the left column of Table 1 then follow from Equation (\ref{eqn:5.1001})
of Remark \ref{rem:5.5}.

On the other hand, let us observe that the Marchenko-Pastur distribution 
of parameter 1 can be written as $\bB_1 ( \, \widetilde{\mu_o} \, )$,
where $\widetilde{\mu_o} = \frac{1}{2} ( \delta_0 + \delta_2 )$. Indeed, 
it is immediately seen by direct calculation that 
$\eta_{\widetilde{\mu_o}} (z) = z/(1-z)$, hence (in view of the above 
Equation (\ref{eqn:2.20})) the measure $\bB_1 ( \, \widetilde{\mu_o} \, )$ 
is determined by the fact that its $R$-transform is $z/(1-z)$; but it is 
well-known that the latter function is exactly the $R$-transform of the 
Marchenko-Pastur (also called free Poisson) distribution of parameter 1 
-- see for instance Section 2.7 of \cite{V00}. The corresponding entry 
in Table 1 is then obtained by just writing that
\[
\phi ( \bB ( \, \widetilde{\mu_o} \, ) ) = 
1+ \phi ( \, \widetilde{\mu_o} \, ) = 1,
\]
where (same as in the preceding paragraph) the fact that  
$\phi ( \, \widetilde{\mu_o} \, ) = 0$ follows from Proposition 
\ref{prop:5.1}, and we used Equation (\ref{eqn:5.1001}) of Remark 
\ref{rem:5.5}.

Finally, it remains to look at the case when $\mu$ is the Cauchy distribution. 
This special $\mu$ has the remarkable property that 
\begin{equation}  \label{eqn:5.5}
\mu^{\boxplus t} = \mu^{\uplus t}, \ \  \forall \, t \geq 0.
\end{equation}
(In order to verify Equation (\ref{eqn:5.5}) one checks that the measures 
on its both sides have the same reciprocal Cauchy transform, which is 
just $F(z) = z+it$.) From (\ref{eqn:5.5}) it follows 
that $\mu$ is fixed by $\bB_t$, for every 
$t \geq 0$, and this in turn makes clear that $\mu \in \bB_t ( \cM )$ for 
all $t \geq 0$ (hence that $\phi ( \mu ) = \infty$, as stated in Table 1). }
\end{remark}

$\ $

$\ $

Serban T. Belinschi: University of Waterloo and IMAR.

Currently at University of Saskatchewan

Address: Department of Mathematics and Statistics, University of Saskatchewan,  

Saskatoon, Saskatchewan, S7N 5E6, Canada.                      

Email: belinschi@math.usask.ca

$\ $

Alexandru Nica: University of Waterloo.

Address: Department of Pure Mathematics, University of Waterloo, 

Waterloo, Ontario N2L 3G1, Canada.                    

Email: anica@math.uwaterloo.ca

\end{document}